\documentclass[english, 12pt,leqno]{amsart}

\usepackage{amsmath}
\usepackage{amsfonts}
\usepackage{amssymb}
\newtheorem{Def}{Definition}

\newtheorem{Cor}[Def]{Corollary}

\newtheorem{Lemma}[Def]{Lemma}
\newtheorem{Prop}[Def]{Proposition}

\newtheorem{Rem}[Def]{Remark}
\newtheorem{Rems}[Def]{Remarks}
\newtheorem{Thm}[Def]{Theorem}
\newtheorem{Item}[Def]{Item}

\title[C$^*$-simple groups]
{C$^*$-simple groups: \\
Amalgamated free products, \\
HNN extensions, \\
and fundamental groups of $3$-manifolds
}

\author{Pierre de la Harpe
and
Jean-Philippe Pr\'eaux
}

\thanks{The authors acknowledge support 
from the Swiss national science foundation}

\address{Pierre de la Harpe~:
Section de math\'ematiques, 
Universit\'e de Gen\`eve, 
C.P.~64, 
CH--1211 Gen\`eve 4. 
}
\email{Pierre.delaHarpe@unige.ch}

\address{Jean-Philippe Pr\'eaux~:
Centre de recherche de l'arm\'ee de l'air,
Ecole de l'air/BA701,
F-13661 Salon-air~; 
\newline
and:
Laboratoire d'analyse topologie et probabilit\'es,
Universit\'e de Provence,
39 rue F. Joliot-Curie,
F-13453 Marseille cedex 13.
}
\email{preaux@cmi.univ-mrs.fr}

\keywords{Groups acting on trees, slender automorphism, 
amalgamated free products, HNN extensions, 
Baumslag-Solitar groups, graphs of groups, 
fundamental groups of $3$-manifolds,
C$\, ^*$-simplicity.}
\subjclass[2000]{22D25, 20E08, 20F65, 57N10}

\date{September 18, 2009, revised June 26, 2011}

\begin{document}

\begin{abstract} 
We establish sufficient conditions for the C$^*$-simplicity
of two classes of groups.
The first class is that of groups acting on trees, 
such as amalgamated free products, HNN-extensions, 
and their non-trivial subnormal subgroups;
for example normal subgroups of Baumslag-Solitar groups.~The 
second class is that of fundamental groups
of compact $3$-manifolds, 
related to the first class by their 
Kneser-Milnor and JSJ-decompo\-sitions.
\par
Much of our analysis deals with conditions 
on an action of a group $\Gamma$ on a tree $T$
which imply  the following three properties:
abundance of hyperbolic elements, 
better called strong hyperbolicity,
minimality, 
both on the tree $T$ and on its boundary $\partial T$,
and faithfulness in a strong sense.
An important step in this analysis is 
to identify automorphism of $T$ which are \emph{slender},
namely such that their fixed-point sets in $\partial T$
are nowhere dense for the shadow topology.
\end{abstract}

\maketitle

\section{\textbf{Introduction
}}
\label{section1}

In the first part of this paper, we analyze actions of groups on trees,
and we  establish that the reduced C$^*$-algebras 
of some of these groups are simple.
In the second part, we apply this to 
fundamental groups of compact $3$-manifolds
and their subnormal subgroups.

Given a group $\Gamma$, recall that its \emph{reduced C$\,{}^*$-algebra}
$C^*_r(\Gamma)$ is the closure for the operator norm
of the group algebra $\mathbf C [\Gamma]$ 
acting by the left-regular representation 
on the Hilbert space $\ell^2(\Gamma)$.
For an introduction to group C$^*$-algebras,
see for example  \cite[Chapter VII]{Davi--96}.
A group is \emph{C$\,{}^*$-simple} if 
it is infinite and if
its reduced C$^*$-algebra 
has no non-trivial two-sided ideals. 
\par

The first examples of groups which have been known to be C$\,{}^*$-simple
are non-abelian free groups \cite{Powe--75}.
Powers' proof applies more generally to the groups defined as follows:
\begin{itemize}
\item[$\bullet$]
a \emph{Powers group}\footnote{Italics indicate a definition,
and the sign $\bullet$ indicate a definition that we wish to emphasize.}
is a group $\Gamma$ which is not reduced to one element
and which has the following combinatorial property: 
\item[]
for any finite subset $F$ in $\Gamma \smallsetminus \{1\}$ 
and for any integer  $n \ge 1$, 
\item[]
there exists a partition $\Gamma = C \sqcup D$ and elements
$\gamma_1,\hdots,\gamma_n$ in $\Gamma$ such that  
\subitem   
$fC \cap C = \emptyset$ for all $f \in F$, and
\subitem
$\gamma_jD \cap \gamma_kD = \emptyset$ 
          for all $j,k \in \{1,\hdots,n\}$, $j \neq k$.
\end{itemize} 
More on Powers groups in Section \ref{sectionPowers},
and in \cite{Harp--07}.
We will also use the following notion:
\begin{itemize}
\item[$\bullet$]
a \emph{strongly Powers group} is a group $\Gamma$ 
such that any subnormal subgroup $N \ne \{1\}$
of $\Gamma$ is a Powers group.
\end{itemize}
Recall that a subgroup $N$ of $\Gamma$ is \emph{subnormal}
if there exists a finite chain of subgroups
$N_0 = N \le N_1 \le \cdots \le N_k = \Gamma$
with $N_{j-1}$ normal in $N_j$ for $j = 1, \hdots, k$.

\begin{Prop}[\textbf{Powers, 1975}]
\label{PowersC*simple}
(i) Powers groups are C$\,{}^*$-simple.
\par
(ii) If $\Gamma$ is a strongly Powers group,
any non-trivial subnormal subgroup of $\Gamma$ is C$\,{}^*$-simple.
\end{Prop}

The proof of Claim (i) is essentially in \cite{Powe--75}; see also \cite{Harp--07};
it is well-known that the converse does \emph{not} hold.
Claim (ii)  follows from (i) and the definition.

The standard scheme to show that a group is a (possibly strongly) Powers group
is to have it act by homeomorphisms on some space $\Omega$
and to call on our Proposition \ref{Powersgeneral}, or a variant thereof.
This applies to the groups of our main result, Theorem \ref{MainResult} below.
Moreover, the proof provides
\emph{informations on the action of $\Gamma$ 
on the appropriate space $\Omega$ which are in our opinion
of independent interest}
(see for example Propositions \ref{PropAmal} and \ref{PropHNN}). 
Before we can state the theorem, 
we need  some notation.

(i)
Let $A,B,C$ be three groups,
let $\iota_A, \iota_B$ be injections of $C$ in $A, B$ respectively,
and let
\begin{equation*}
\Gamma \, = \, A \ast_C B \, = \, 
\langle A,B \hskip.1cm \vert \hskip.1cm \iota_A(c) = \iota_B(c)
\hskip.2cm \forall c \in C \rangle
\end{equation*}
be the corresponding amalgam.
We identify $C$ to a subgroup of $A$ and of $B$;
we define inductively a decreasing sequence 
$C_0 \supset C_1 \supset \cdots$ of subgroups of $C$
by $C_0 = C$ and 
\begin{equation*}
C_{k+1} \, = \, 
\big( \bigcap_{a \in A} a^{-1}C_k a \big)  
\cap
\big( \bigcap_{b \in B} b^{-1}C_k b \big)
\end{equation*}
for $k \ge 0$.
For the geometrical meaning of the condition $C_k = \{1\}$,
see Subsection \ref{ssamalgamated}.
\par

(ii)
Let $G$ be a group, let $\theta$ be an isomorphism 
from a subgroup $H$ of $G$ to some subgroup of $G$,
and let
\begin{equation*}
\Gamma \, = \,  HNN(G, H, \theta)
\, = \,
\langle G, \tau \hskip.1cm \vert \hskip.1cm
\tau^{-1}h\tau = \theta(h) \hskip.2cm
\forall h \in H \rangle 
\end{equation*}
be the corresponding HNN-extension.
We define inductively a decreasing sequence
$H_0 \supset H_1 \supset \cdots$ of subgroups of $H$
by $H_0 = H$ and
\begin{equation*}
\aligned
H'_k \, &= \, H_k \cap \tau^{-1} H_k \tau \, = \, 
H_k \cap \theta( H_k )
\\
H_{k+1} \, &= \, 
\Big( 
    \bigcap_{g \in G} g  H'_k  g^{-1} 
\Big)  
\cap
\tau \Big( 
    \bigcap_{g \in G} g  H'_k  g^{-1} 
\Big)   \tau^{-1}
\endaligned
\end{equation*}
for $k \ge 0$.
For the geometrical meaning of the condition $H_k = \{1\}$,
see Subsection \ref{ssHNN}. 
\par

(iii)
For  $m,n \in \mathbf Z$,
the corresponding \emph{Baumslag-Solitar group} is 
\begin{equation*}
\aligned
BS(m,n) \, &= \, 
\langle \tau,b \hskip.1cm \vert \hskip.1cm \tau^{-1}b^m\tau = b^n \rangle
\\
\, &= \,
HNN(b^{\mathbf Z}, b^{m\mathbf Z}, \hskip.1cm
b^{mk} \longmapsto b^{nk}) .
\endaligned
\end{equation*}
(In this case, $H_k \ne \{1\}$ for any $k \ge 0$;
see the proof of Claim (iii) of Theorem \ref{MainResult},
at the end of Subsection \ref{ssHNN}.)
\par

(iv) Let $M$ be a $3$-manifold.
Denote by $\widehat M$ the manifold obtained from $M$ by filling all $2$-spheres
in $\partial M$ (if any) with $3$-balls,
and by $\Gamma$ the fundamental group $\pi_1(M) \approx \pi_1(\widehat M)$.
For our phrasing of the last claim in the next theorem,
we rely on Perelman's proof of Thurston's geometrisation conjecture;
see Section \ref{sectionparticularcases} for some relevant definitions.

\begin{Thm} 
\label{MainResult}
(i) 
For a countable amalgam
$\Gamma \, = \, A \ast_C B$ to be a strongly Powers group,
it is sufficient that
$[A : C] \ge 3$, $[B : C] \ge 2$, 
and $C_k \, = \, \{1\}$ for some $k \ge 0$.
\par

(ii) 
For a countable HNN-extension 
$\Gamma =  HNN(G, H, \theta)$ to be a strongly Powers group,
it is sufficent that
$H \subsetneqq G$, $\theta(H) \subsetneqq G$,
and $H_k = \{1\}$ for some $k \ge 0$.
\par

(iii) 
A Baumslag-Solitar group $B(m,n)$ 
is a strongly Powers group 
if and only if  it is C$\,^*$-simple,
if and only if 
$\min\{\vert m\vert , \vert n \vert\} \ge 2$
and $\vert m \vert \ne \vert n \vert$.
\par

(iv)
Let $M$ be a $3$-manifold 
which is compact, connected, and orientable,
and let $\Gamma = \pi_1(M) \approx \pi_1(\widehat M)$.
If $\widehat M$ is a Seifert manifold, 
or if the interior of $M$ is a Sol-manifold,
then $\Gamma$ is not C$\,{}^*$-simple.
Otherwise, $\Gamma$ is a strongly Powers group.
\par
In particular, let $\Gamma$ be the group of a tame knot $K$ in the $3$-sphere.
Then $\Gamma$ is a strongly Powers group
if and only if $\Gamma$ is C$\,{}^*$-simple,
if and only if $K$ is not a torus knot.
\end{Thm}

A well-known particular case of Claim (i) is that of free products  \cite{PaSa--79}.
We agree that
\begin{itemize}
\item[$\bullet$]
a free product $A \ast B$ is \emph{non-trivial}
if neither $A$ nor $B$ is the group with one element.
\end{itemize}
Recall that the 
\begin{itemize}
\item[$\bullet$]
\emph{infinite dihedral group}
is the free product of two groups of order two,
\end{itemize}
and that it is the only amenable non-trivial free product.

\begin{Cor}
\label{Freeproducts}
For a non-trivial free product $\Gamma = A \ast B$, 
the three following conditions are equivalent:
$\Gamma$ is not isomorphic to the infinite dihedral group,
$\Gamma$ is C$\, ^*$-simple,
$\Gamma$ is a strongly Powers group.
\end{Cor}

At first sight, Corollary \ref{Freeproducts} 
follows from Claim (i) for countable groups only.
But a group $\Gamma$ is C$^*$-simple 
[respectively Powers, strongly Powers] as soon as,
for any countable subgroup $\Gamma_0$ of $\Gamma$,
there exists a countable C$^*$-simple 
[respectively Powers, strongly Powers] subgroup
of $\Gamma$ containing $\Gamma_0$;
see \cite[Proposition 10]{BeHa--00} and \cite[Lemma 3.1]{BoNi--88}.
It follows that Corollary~\ref{Freeproducts}
holds as stated above.
\par

To emphasize the benefit of proving the \emph{strong} Powers property,
we state the following corollary of Claim (iii):

\begin{Cor}
\label{SBS}
For $BS(m,n)$ as in Theorem~\ref{MainResult},
let $SBS(m,n)$ denote the kernel of 
the homomorphism $\sigma : BS(m,n) \longrightarrow \mathbf Z$
defined by $\sigma(b) = 0$ and $\sigma(\tau) = 1$.
\par
If $\min\{\vert m \vert, \vert n \vert\} \ge 2$
and $\vert m \vert \ne \vert n \vert$,
then $SBS(m,n)$ is a strongly Powers group.
\end{Cor}

A \textbf{first part of our paper} contains the proof of the three first claims
of Theorem~\ref{MainResult}.
Much of it is a reminder, in a form adapted to our purpose:
on Powers group in 
Section \ref{sectionPowers},
on groups acting on trees
in Section \ref{sectiongroupsontrees},
and on the fundamental group $\Gamma$
of a graph of groups $\mathbf G  = (G,Y)$,
which acts on its Bass-Serre tree $T$, in
Section \ref{sectionBS}.
We recall known criteria for $T$ to be ``small''
(a vertex, or a linear tree),
and for the action to be minimal.
Most of this can be found in papers by H.~Bass.
In Section~\ref{section5}, we analyze further
the two standard examples
and prove  (i) to (iii) of Theorem~\ref{MainResult};
though it does not  follow from Claim (ii),
Claim (iii) can be shown by similar arguments.
\par

Let us now comment on the improvements of these
over previously known results.
It is for these improvements that we have found it useful
to write up our reminders in great detail.
\par

Claim (i) of Theorem \ref{MainResult}
is an improvement in two ways.
(a)~It establishes a property of any subnormal subgroup $N \ne \{1\}$, 
and not only of $N = \Gamma$.
(b)~Its hypothesis are stated in terms of $A$, $B$, and $C$ only;
on the contrary, 
with the weaker conclusion ``$\Gamma$ is a Powers group'',
Proposition 10 in \cite{Harp--85} and Corollary 4.6 in \cite{Ivan}
have hypothesis stated in terms
of the action of $\Gamma$ on the edge set of its Bass-Serre tree
(more precisely, 
this action should be strongly faithful, 
as defined in Section~\ref{sectionPowers}).
Our conditions 
are also weaker than those of \cite{Bedo--84}.
\par

Similar remarks hold for Claim (ii) and previous results,
such as Proposition 11 in \cite{Harp--85}.
The C$\,{}^*$-simplicity of $B(m,n)$ for $m,n$ as in Claim (iii)
is due to Nikolay Ivanov \cite[Theorem 4.9]{Ivan}.
As much as we see, Corollary \ref{SBS} does not follow from Ivanov's arguments. 

\vskip.2cm

In a \textbf{second part of our paper},
we consider the fundamental group $\Gamma$
of a  connected compact $3$-manifold $M$,
and we prove Claim (iv) of Theorem~\ref{MainResult}.
Whenever $\Gamma$ is the fundamental group of the graph of groups
given by some Kneser-Milnor or JSJ decomposition of $M$,
the group $\Gamma$ has a canonical action
on the corresponding Bass-Serre tree, and the results of Part I apply.
\par

In this introduction, 
rather than describing  precisely the case of $3$-manifold groups, 
let us first record the much simpler situation of fundamental groups of surfaces.
Let us recall that a group is 
\begin{itemize}
\item[$\bullet$]
\emph{icc} if it is infinite and if all its conjugacy classes
other than $\{1\}$ are infinite.
\end{itemize}
Let $S$ be a connected compact surface, with fundamental group $\Gamma$.
It follows from the classification of surfaces and from elementary arguments
(compare with \cite{HaPr--07}) that the following three properties are equivalent:
\begin{itemize}
\item[(i)]
$S$ is not diffeomorphic to a disc, a sphere, a projective plane, 
an annulus, a M\"obius band, a $2$-torus, or a Klein bottle; 
\item[(ii)]
$\Gamma$ is icc,
or equivalently the von Neumann algebra of $\Gamma$
is a factor of type $II_1$;
\item[(iii)]
$\Gamma$ is C$^*$-simple.
\end{itemize}
If $\Gamma$ is now the fundamental group 
of a connected compact $3$-manifold,
the equivalence of (ii) and (iii) does not carry over.
Indeed, the fundamental group of a Sol-manifold can be icc,
and is never C$^*$-simple.
In some sense, Claim (iv) of Theorem \ref{MainResult} 
provides the $3$-dimensional analogue of these equivalences.

\vskip.2cm

We are grateful to Luc Guyot for suggesting 
that most of our propositions establish not only that 
some $\Gamma$ is a Powers group,
but also that any subnormal subgroup $N \ne \{1\}$ 
is a Powers group,
and for many other comments
on  preliminary versions of this paper.
We are also grateful to
Laurent Bartholdi, Bachir Bekka, Ken Dykema,  
Yves Stalder, and Nicolas Monod 
for  helpful remarks.

\section{\textbf{A reminder on Powers groups
}}
\label{sectionPowers}

The main result of this section is Proposition \ref{Powersgeneral}, 
which is the standard tool for showing groups to be Powers groups.
Before this, we define the relevant notions and establish two lemmas.

An action of a group $\Gamma$ on a set $\Omega$ is
\begin{itemize}
\item[$\bullet$]
\emph{faithful} if, 
for any $\gamma \in \Gamma$, $\gamma \ne 1$, 
there exists $\omega \in \Omega$ such that
$\gamma(\omega) \ne \omega$, and
\item[$\bullet$]
\emph{strongly faithful} if, 
for any finite subset $F$ of $\Gamma$ not containing $1$,
there exists $\omega \in \Omega$ such that 
$\gamma (\omega) \ne \omega$ for all $\gamma \in F$.
\end{itemize}
Let  $\Omega$ be a Hausdorff space 
and let the action of $\Gamma$ be by homeomorphisms;
if the action is strongly faithful 
and if $F$, $\omega$ are as above,
observe that there exists 
a neighbourhood $V$ of $\omega$ in $\Omega$
such that $\gamma (V) \cap V = \emptyset$ 
for all $\gamma \in F$.
The proof of the following lemma is straightforward.

\begin{Lemma}
\label{defsf}
Let $\Gamma$ be a group acting by homeomorphisms
on a topological space $\Omega$
and let $L$ be a $\Gamma$-invariant dense subspace of $\Omega$.
\par

The action of $\Gamma$ on $L$ is strongly faithful
if and only if the action of $\Gamma$ on $\Omega$ is strongly faithful.
\end{Lemma}

An action of a group $\Gamma$ by homeomorhisms
on a topological space $\Omega$ is 
\begin{itemize}
\item[$\bullet$]
\emph{minimal} if all the $\Gamma$-orbits are dense in $\Omega$.
\end{itemize}
A homeomorphism $\gamma$ of a Hausdorff space $\Omega$ is
\begin{itemize}
\item[$\bullet$]
\emph{hyperbolic} if it has two fixed points  $\alpha, \omega \in \Omega$ such that,
for any neighbourhoods $U$ of $\alpha$ and $V$ of $\omega$,
and for $n$ large enough,
$\gamma^n(\Omega \smallsetminus U) \subset V$
and $\gamma^{-n}(\Omega \smallsetminus V) \subset U$.
\end{itemize}
An action of a group $\Gamma$ on $\Omega$ is
\begin{itemize}
\item[$\bullet$]
\emph{strongly hyperbolic}
if $\Gamma$ contains two hyperbolic homeomorphisms
which are \emph{transverse}, namely
without common fixed point.
\end{itemize}
If the action is strongly hyperbolic, 
there is an abundance of transverse pairs:

\begin{Lemma}
\label{defsh}
Let $\Gamma$ be a group acting strongly hyperbolically 
on a topological space $\Omega$.
\par

For any hyperbolic element $\gamma_0 \in \Gamma$,
there exist infinitely many hyperbolic elements $\gamma_n$, $n \ge 1$,
such that the $\gamma_n$, $n \ge 0$, are pairwise transverse.
\end{Lemma}

\noindent 
\emph{Proof.}~Let $\alpha, \beta \in \Gamma$ be hyperbolic and transverse.
Upon replacing $\alpha$ with $\beta^k \alpha \beta^{-k}$ for an appropriate $k$,
we can assume furthermore that $\gamma_0$ and $\alpha$ are transverse.
\par
Set $\delta_m = \alpha^m \gamma_0 \alpha^{-m}$ for all $m \ge 1$.
Then there exists a subsequence $(\gamma_n)_{n \ge 1}$
of  $(\delta_m)_{m \ge 1}$ such that the $\gamma_n$, $n \ge 0$,
are pairwise transverse.
\hfill $\square$

\vskip.2cm

If $\Gamma$ acts on $\Omega$, we denote by $L_{\Gamma}$ 
the set of points $\omega \in \Omega$ such that $\omega = s(\gamma)$
for some hyperbolic element $\gamma \in \Gamma$
(note that ``L'' is as in ``Limit set'').
Observe that $L_{\Gamma}$ is $\Gamma$-invariant,
and infinite in case the action is strongly hyperbolic.

In the following proposition,
Claim (ii) is  a reformulation
of Proposition~11 in \cite{Harp--07}.

\begin{Prop}
\label{Powersgeneral}
Let $\Gamma$ be a group which acts by homeomorphisms
on a Hausdorff topological space $\Omega$.
Assume that the action
is strongly hyperbolic  on $\Omega$
and strongly faithful on $L_{\Gamma}$.
Then:
\begin{itemize}
\item[(i)]
any non-empty $\Gamma$-invariant closed subset of $\Omega$
contains $\overline{L_{\Gamma}}$;
\item[(ii)]
$\Gamma$ is a Powers group.
\end{itemize}
More generally:
\begin{itemize}
\item[(iii)]
any subnormal subgroup $N \ne \{1\}$ in $\Gamma$
which contains a hyperbolic element is a Powers group.
\end{itemize}
\end{Prop}

\noindent 
\emph{Proof.}~
(i)
Let $C$ be a non-empty  $\Gamma$-invariant closed
subset of $\Omega$, and choose $\xi \in C$.
Let $\eta \in L_{\Gamma}$; 
choose a hyperbolic element $\gamma \in \Gamma$ 
such that $\eta = \omega(\gamma)$.
\par
Assume first that  $\xi \ne \alpha(\gamma)$.
We have $\eta = \lim_{n \to \infty} \gamma^n (\xi)$;
since $\gamma^n(\xi) = \alpha(\gamma^n \alpha \gamma^{-n})$,
this implies $\eta \in C$.
Assume next that $\xi = \alpha(\gamma)$.
By Lemma \ref{defsh}, 
there exists a hyperbolic element $\gamma' \in \Gamma$
transverse to $\gamma$.
We have first $\omega(\gamma') \in C$ as in the previous case,
and then  
$\eta = \lim_{n \to \infty} \gamma^n \big( \omega(\gamma') \big) \in C$.
\par

This implies $C \supset L_{\Gamma}$, and consequently
$C = \overline{L_{\Gamma}} = \partial T$.

\vskip.2cm

(ii)
Let $F \subset \Gamma \smallsetminus \{1\}$
be a finite subset and $n \ge 1$ be an integer,
as in the definition of ``Powers group''.
\par

By hypothesis, there exists a hyperbolic element $\gamma_0$
in $\Gamma$ and a neighbourhood $C_{\Omega}$
of $\omega(\gamma_0)$ in $\Omega$ 
such that $f(C_{\Omega}) \cap C_{\Omega} = \emptyset$
for all $f \in F$.
By the proof of Lemma \ref{defsh}, 
there exist $\gamma_1, \hdots, \gamma_n \in \Gamma$
such that $\gamma_0, \gamma_1, \hdots, \gamma_n$
are pairwise transverse hyperbolic elements.
Choose neighbourhoods $A_j$ of $\alpha(\gamma_j)$
and $\Omega_j$ of $\omega(\gamma_j)$, $j = 1,\hdots,n$,
such that $A_1, \Omega_1, \hdots, A_n, \Omega_n$
are pairwise disjoint.
Upon replacing $\gamma_j$ by a large enough power of itself,
we can assume that 
$\gamma_j^k(\Omega \smallsetminus A_j) \subset \Omega_j$ and
$\gamma_j^{-k}(\Omega \smallsetminus \Omega_j) \subset A_j$
for all $k \ge 1$.
Upon replacing now $\gamma_j$ by 
$\gamma_0^{\ell} \gamma_j \gamma_0^{-\ell}$ for $\ell$ large enough,
we can assume furthermore that
$A_1 \cup \Omega_1 \cup \cdots \cup A_n \cup \Omega_n \subset C_{\Omega}$.
\par

Choose $\omega_0 \in \Omega$. Set
$C = \{ \gamma \in \Gamma \hskip.1cm \vert \hskip.1cm 
\gamma(\omega_0) \in C_{\Omega} \}$
and $D = \Gamma \smallsetminus C$.
For $f \in F$, we have $fC \cap C = \emptyset$
since $f(C_{\Omega}) \cap C_{\Omega} = \emptyset$.
For $j \ne k$ in $\{1,\hdots,n\}$, 
we have $\gamma_j D \cap \gamma_k D = \emptyset$
since $\Omega_j \cap \Omega_k = \emptyset$.
This ends the proof of Claim (ii).

\vskip.2cm

(iii) 
Consider first the particular case of a group $N$ which is normal in $\Gamma$,
so that $L_N$ is a $\Gamma$-invariant subset of $L_{\Gamma}$.
By hypothesis, there exists a hyperbolic element $\gamma_1 \in N$,
say with source $\alpha_1$ and sink $\omega_1$;
thus $L_N \ne \emptyset$, and therefore $\overline{L_N} = \overline{L_{\Gamma}}$.
We will now show that the action of $N$ on $\Omega$ 
(equivallently on $\overline{L_N}$) is strongly hyperbolic.
\par

Since $L_N$ is infinite, 
there exists a hyperbolic element $\gamma_2 \in N$,
say with source $\alpha_2$ and sink $\omega_2$,
such that $\{\alpha_2, \omega_2\} \ne \{\alpha_1, \omega_1\}$.
Upon replacing $\gamma_2$ by $\gamma_2^{-1}$,
we may assume that 
$\alpha_2 \notin \{\alpha_1,\omega_1\}$.
If $\omega_2 \notin \{\alpha_1,\omega_1\}$,
there is nothing left to prove; 
upon replacing $\gamma_1$ by $\gamma_1^{-1}$ if necessary,
we can assume from now that $\omega_2 = \omega_1$,
and we denote this point by $\omega$.
There cannot be any point in $L_N$ fixed by $N$.
Hence there exists $\delta \in N$ 
such that $\delta(\omega) \ne \omega$.
Upon exchanging $\gamma_1$ and $\gamma_2$,
we can assume that $\delta(\alpha_1) \ne \omega$.
\par

We claim that all the orbits of $N$ on $\Omega$ are infinite.
Since any orbit of $\gamma_1^{\mathbf Z}$ 
on $\Omega \smallsetminus \{\alpha_1,\omega\}$ is infinite,
and similarly for $\gamma_2^{\mathbf Z}$
on $\partial T \smallsetminus \{\alpha_2,\omega\}$, 
the only point we have to check is that 
the orbit $N(\omega)$ is infinite.
Observe that $\gamma_3 \Doteq \delta \gamma_1 \delta^{-1}$
is hyperbolic, that $\alpha(\gamma_3) = \delta(\alpha_1)$,
and that $\omega(\gamma_3) = \delta(\omega)$.
Since $\omega \notin \{\alpha(\gamma_3),\omega(\gamma_3)\}$,
the $\gamma_3^{\mathbf Z}$-orbit of $\omega$ is infinite.
\emph{A fortiori}, the $N$-orbit of $\omega$ is infinite.
\par

It is a general fact that, if a group $N$
acts on a set $\Omega$ in such a way that all its orbits are infinite,
and if $F$ is any finite subset of $\Omega$,
there exists $\gamma \in N$ such that $F$ and $\gamma(F)$
are disjoint; see\footnote{
This is a straightforward consequence of the following lemma
of B.H.~Neumann \cite{NeBH--54}:
a group cannot be covered by finitely many cosets
with respect to subgroups of infinite index.
Other and ``new''  proof: Lemma 2.4 in \cite{PoVa--08}.
} 
Lemma 2.3 in \cite{NePM--76}.
In our case, it implies that we can choose $\gamma \in N$ such that
$\gamma(\{\alpha_1,\omega\}) \cap \{\alpha_1,\omega\} = \emptyset$,
namely such that the hyperbolic elements 
$\gamma_1$ and $\gamma \gamma_1 \gamma^{-1}$ are transverse.
\par

Note that the action of $N$ on $\overline{L_N}$ is strongly faithful
(see Lemma \ref{defsf}).
The proof of Claim (ii) shows now that $N$ is a Powers group.
\par

Claim (iii) for the general case ($N \ne \{1\}$ subnormal)
follows from the proof in the special case ($N \ne \{1\}$ normal).
\hfill $\square$

\vskip.2cm

Concerning Powers groups,
let us furthermore quote some general facts;
for these and for complements, we refer to \cite{Harp--07}.
\par

\begin{Item}  
\label{ItemPowers}
The following facts about Powers groups and C$^*$-simple groups are well-known.
\begin{itemize}
\item[(i)]
A Powers group is C$^*$-simple, 
and moreover its reduced C$^*$-algebra has a unique trace 
(see Proposition \ref{PowersC*simple}). 
\item[(ii)]
Powers groups are icc.
\item[(iii)]
Powers group have non-abelian free subgroups \cite{BrPi}.
\item[(iv)]
Torsion-free Gromov-hyperbolic groups are Powers groups;
in particular, non-abelian free groups are Powers groups.
\item[(v)]
There are uncountably many countable groups $\Gamma$
with pairwise non-isomorphic simple reduced C$^*$-algebras
(this is Corollary~9 in \cite{AkLe--80}, building up on \cite{McDu--69}).
\item[(vi)]
The only amenable normal subgroup of a Powers group
(or more generally of a C$^*$-simple group)
is the trivial group $\{1\}$.
\item[(vii)]
There are  C$^*$-simple groups which are not Powers groups
(for example any
direct product of non-abelian free groups).
\end{itemize}
\end{Item}

We do not know any example of a Powers group 
which is not a strongly Powers group, as defined in the introduction;
similarly, the following question is open for us:
\begin{itemize}
\item[(Q)]
\emph{does there exist a pair $(\Gamma,N)$
of a C$\,^*$-simple group $\Gamma$ 
and a non-trivial normal subgroup $N$ 
which is not C$\, ^*$-simple?}
\end{itemize}
We are grateful to Bachir Bekka who has observed to us that,
as a consequence of \cite{Pozn}, if such a pair $(\Gamma,N)$ exists,
then $N$ cannot be linear (\emph{a fortiori} $\Gamma$ cannot be linear).
Indeed, suppose that a C*-simple  group $\Gamma$
would contain a linear non-C*-simple non-trivial normal subgroup $N$; 
then, by Poznansky's results, the amenable radical $R$ of $N$ would be non-trivial;
but amenable radicals are characteristic subgroups, 
so that $R$ would be normal in $\Gamma$,
and this is impossible by Item \ref{ItemPowers}.vi. 
The question is closely related to that of
asking wether there exists a group $\Gamma$ 
which is \emph{not} C$^*$-simple 
and which does not contain any non trivial amenable normal subgroup
\cite{BeHa--00}.

\section{\textbf{On groups acting on trees and their boundaries
}}
\label{sectiongroupsontrees}

Let $X$ be a \emph{graph}, with vertex set $V(X)$ and edge set $E(X)$.
We follow Bass \cite{Bass--93} and Serre \cite{Serr--77};
in particular, each geometric edge of $X$
corresponds to a pair $\{e,\overline{e}\} \in E(X)$,
with $\overline{e} \ne e$ and $\overline{\overline{e}} = e$.
Each edge $e \in E(X)$ has a \emph{source} $s(e) \in V(X)$ and a 
\emph{terminus} $t(e) \in V(X)$.
We denote by $d(x,y)$ the \emph{combinatorial distance}
between two vertices $x,y \in V(X)$.
\par

In case of a \emph{tree}, we write $T$ rather than $X$,
and we agree that $V(T) \ne \emptyset$.
Recall that a \emph{ray} in $T$ is a subtree 
with vertex set $(x_n)_{n \in \mathbf N}$, such that
\begin{equation*}
d(x_m,x_n) \, = \,  \vert m-n \vert \hskip.5cm
\text{for all} \hskip.2cm m,n \in \mathbf N.
\end{equation*}
The enumeration of the vertices of a ray will \emph{always}
be such that the above condition holds.
Two rays with vertex sets 
$(x_n)_{n \in \mathbf N}$ and $(y_n)_{n \in \mathbf N}$
are \emph{cofinal} if there exists some $k \in \mathbf Z$
such that $y_n = x_{n+k}$ for all $n$ large enough.
\begin{itemize}
\item[$\bullet$]
The \emph{boundary} $\partial T$ of a tree $T$
is the set of cofinal classes of rays in $T$.
\end{itemize}
Observe that, for $\partial T$ to be non-empty,
it suffices that $T$ does not have any  vertex of degree $\le 1$
(since we have agreed that $V(T) \ne \emptyset$).
\par

The set $\partial T$ has a natural topology defined as follows.
For any edge $e \in E(T)$, 
the \emph{shadow} $(\partial T)_e$ of $e$ in $\partial T$
is the subset of $\partial T$
represented by rays $(x_n)_{n \in \mathbf N}$
such that $d(t(e),x_n) < d(s(e),x_n)$ for all $n$ large enough.
The family of shadows $\left\{ (\partial T)_e \right\}_{e \in E(T)}$
generates a topology on $\partial T$
which is Hausdorff and totally disconnected;
if $T$ is countable, this topology is moreover metrisable.
Given any  $x_0 \in V(T)$,
this topology coincides\footnote{Similarly,
one can define two topologies on the disjoint union
$\overline{T} = V(T) \sqcup \partial T$,
one using appropriate shadows and 
the other by identifying $\overline{T}$ 
to the inverse limit of the balls
$\{x \in V(T) \hskip.1cm \vert \hskip.1cm d(x_0,x) \le n\}$.
These two topologies coincide if and only if
the tree $T$ is locally finite,
but their restrictions on $\partial T$ coincide in all cases.
The shadow topology makes $\overline{T}$ a compact space,
but $\partial T$ need not be closed,
and therefore need not be compact
(unless $T$ is locally finite). More on this in \cite{MoSh--04}.
}
with the inverse limit topology
on the set $\partial T$,  
identified with the inverse limit of the discrete spaces, or \emph{spheres},
$S(x_0,n) = \{x \in V(T) \hskip.1cm \vert \hskip.1cm d(x_0,x) = n\}$.
It is known that $\partial T$ is compact if and only if $T$ is locally finite
(and we do consider below trees \emph{which are not} locally finite).
See the last exercise in Section I.2.2 of \cite{Serr--77},
Section I.8.27 in \cite{BrHa--99},
and Section 4 in \cite{MoSh--04}.

\begin{itemize}
\item[$\bullet$]
A \emph{pending ray} in a tree $T$ is a ray with vertex set
$\left(x_n\right)_{n \ge 1}$ such that  $x_n$ has degree $2$ in $T$
for all $n$ large enough.
\end{itemize}
A ray cofinal with a pending ray is also pending.
The classes of the pending rays are precisely 
the isolated points in $\partial T$;
in particular, the topological space $\partial T$
is \emph{perfect} if the tree $T$ does not have any pending ray.

\begin{Prop}
\label{Baire}
For any countable tree $T$, the boundary $\partial T$ is a Baire space.
\end{Prop}

\emph{Proof.}
Choose a vertex $x_0 \in V(T)$.
Since $T$ is countable, 
each of the spheres $S(x_0,n)$ defined above is countable,
and therefore is a polish space for the discrete topology.
Since the inverse limit of a countable inverse system of polish spaces 
is polish \cite[$\S$ 6]{Bou--TG9},
$\partial T$ is polish
(``countable'' is important, because a product 
of uncountably many polish spaces is not polish in general).
Since polish spaces are  Baire spaces 
\cite[$\S$~5]{Bou--TG9},
this ends the proof.
\hfill $\square$

\vskip.2cm

An action of a group on a tree $T$ is faithful, or strongly faithful,
if it is the case for the action on the vertex set $V(T)$;
see Section \ref{sectionPowers}.
\par

Any automorphism $\gamma$ of a tree $T$ induces
a homeomorphism of the boundary $\partial T$, again denoted by $\gamma$.
An automorphism $\gamma$ of a tree $T$ is 
\begin{itemize}
\item[$\bullet$]
\emph{slender} if $\partial T \ne \emptyset$\footnote{As 
our trees are assumed to be non-empty,
observe  that a tree on which a group acts 
without globally fixed vertex and without globally fixed geometric edge
has infinite diameter, and more precisely non-empty boundary.
}
and if the fixed point set of $\gamma$ in $\partial T$ has empty interior.
\end{itemize}
The action of $\Gamma$ on $T$ is \emph{slender}
if any $\gamma \in \Gamma$, $\gamma \ne 1$, is slender.
Observe that a slender action is faithful.

\begin{Cor}
\label{slender}
Let $\Gamma$ be a  group which acts by automorphisms 
on a countable tree $T$.
\par
If the action on $T$ is slender, 
then the action on $\partial T$ is strongly faithful.
\end{Cor}

\noindent 
\emph{Proof.}~
Let $F$ be a finite subset of $\Gamma \smallsetminus \{1\}$,
as in the definition of ``strongly faithful''.
Since the action is slender,
the fixed point set $(\partial T)^\gamma$ has no interior 
for all $\gamma \in F$, 
so that there exists by Proposition \ref{Baire}
a point $\xi \in \partial T$
such that $\gamma(\xi) \ne \xi$ for all $\gamma \in F$.

[Observe that
$\partial T \smallsetminus 
\bigcup_{\gamma \in \Gamma \smallsetminus \{1\}} (\partial T)^{\gamma}$
is dense, and therefore non-empty, 
so that there exists $\xi \in \partial T$
such that $\gamma(\xi) \ne \xi$ 
for all $\gamma \in \Gamma \smallsetminus \{1\}$.
We shall not use this below.]
\hfill $\square$

\vskip.2cm

Since \cite{Tits--70}, we know that
an automorphism $\gamma$ of a tree $T$
can be of three different kinds: 
\begin{itemize}
\item[$\bullet$]
It is an \emph{inversion} if there exists $e \in E(T)$
such that $\gamma(e) = \overline{e}$.
\item[$\bullet$]
It is \emph{elliptic} if it has at least one fixed point in $V(T)$.
\item[$\bullet$]
It is \emph{hyperbolic} if it is not an inversion and if 
\begin{equation*}
d(\gamma) \Doteq \min\{d(y,\gamma(y)) 
\hskip.1cm \vert \hskip.1cm y \in V(T) \} \ge 1.
\end{equation*}
If $\gamma$ is hyperbolic, it has an \emph{axis},
which is a linear subtree $T_{\gamma}$, such that  
$x \in V(T_{\gamma})$ if and only if
$d(x,\gamma(x)) = d(\gamma)$,
and $d(\gamma)$ is called the \emph{translation length} of $\gamma$.
\end{itemize}
Note that a tree $T$ with $\partial T = \emptyset$ 
does not have any hyperbolic automorphism.
When $\partial T \ne \emptyset$ note that an automorphism of $T$
is hyperbolic if and only if the induced homeomorphism of $\partial T$
is hyperbolic in the sense of Section \ref{sectionPowers}.

We translate here by 
\begin{itemize}
\item[$\bullet$]
\emph{linear subtree} 
\end{itemize}
what is a ``cha\^{\i}ne doublement infinie'' in \cite{Tits--70}
and a ``droit chemin'' in \cite{Serr--77},
namely a subtree of which the vertex set
is of the form $(x_n)_{n \in \mathbf Z}$,
with $d(x_m,x_n) = \vert m - n \vert$ for all $m,n \in \mathbf Z$.
Linear trees are defined accordingly.

\begin{Rems}
\label{Remarks}
Let $T$ be a tree. \par
(i)
An automorphism $\gamma$ of $T$ which is hyperbolic
has exactly two fixed points on $\partial T$
which are its \emph{source} $\alpha(\gamma)$
and its \emph{sink} $\omega(\gamma)$,
and the infinite cyclic group $\gamma^{\mathbf{Z}}$
acts freely on 
$\partial T \smallsetminus \{\alpha(\gamma), \omega(\gamma) \}$.
If $\gamma$ is a hyperbolic automorphism of a tree
of which the boundary consists of more than two points
(and therefore of which the boundary is infinite), 
then $\gamma$ is always slender, 
since its fixed point set in $\partial T$ 
consists of two non-isolated points.
\par
(ii)
If $\partial T \ne \emptyset$
and if $\gamma$ is an elliptic automorphism of $T$
with fixed point set $V(T)^{\gamma}$ of finite diameter,
then $(\partial T)^{\gamma} = \emptyset$,
and in particular $\gamma$ is slender.
\par
(iii)
Let $\gamma$ be an elliptic automorphism of  $T$
such that $(\partial T)^\gamma \ne \emptyset$.
Consider a ray with vertex set $(x_n)_{n \ge 0}$
of which the origin $x_0$ is fixed by $\gamma$
and which represents a boundary point fixed by $\gamma$;
then each of the vertices $x_n$ is fixed by $\gamma$.
\par
(iv)
If $T$ is the Bass-Serre tree corresponding 
to the Baumslag-Solitar group $\Gamma = BS(m,n)$ 
for some $m,n \in \mathbf Z$
(as in Subsection \ref{ssHNN} below), 
there exist  elliptic automorphisms $\gamma \in \Gamma$ of $T$
with fixed-point sets of infinite diameter
\cite[\emph{Exemple} 4.2]{Stal--06}.
\par
(v)
Elliptic automorphisms need not be slender.
\end{Rems}

For \emph{(v)}, consider for example the regular tree $T$ of degree $3$,
a vertex $x_0$ of $T$, 
and the three isomorphic connected components $T_1$, $T_2$, $T_3$
obtained from $T$ by deleting $x_0$ and the incident edges.
An appropriate transposition $\sigma$ of $T_1$ and $T_2$ 
fixing $T_3$ is not slender,
since $(\partial T)^{\sigma} = \partial T_3$ is open non-empty, 
but a cyclic permutation $\gamma$ of $T_1$, $T_2$, $T_3$
of order~$3$ is slender since $(\partial T)^{\gamma} = \emptyset$.

\vskip.2cm

Two hyperbolic automorphisms of a tree are
\begin{itemize}
\item[$\bullet$]
\emph{transverse} if they do not have any common fixed point in $\partial T$
(see Section \ref{sectionPowers}),
equivalently if the intersections of their axis is finite (possibly empty).
\end{itemize}
An action of a group $\Gamma$ on a tree $T$ is
\begin{itemize}
\item[$\bullet$]
\emph{strongly hyperbolic} if
$\Gamma$ contains a pair of transverse hyperbolic elements,
namely if and only if the corresponding action of $\Gamma$
on the topological space $\partial T$ is strongly hyperbolic
in the sense of Section \ref{sectionPowers}.
\end{itemize}
Lemma \ref{defsh} shows that, 
if the action of $\Gamma$ on $T$ is strongly hyperbolic,
then $\Gamma$ contains an abundance of transverse pairs of hyperbolic elements.

There is a notion of an \emph{irreducible} action of a group
on a $\Lambda$-tree, where $\Lambda$ is an ordered abelian group;
in the special case of a $\mathbf{Z}$-tree, these irreducible actions
are our strongly hyperbolic actions, as it follows
from Proposition 3.7 of Chapter 3 in \cite{Chis--01}.
\par

An action of a group $\Gamma$ on  a tree $T$ is
\begin{itemize}
\item[$\bullet$] 
\emph{minimal}
if there does not exist any proper $\Gamma$-invariant subtree in $T$.
\end{itemize}
Since the same word ``minimal'' is used for actions on trees as above
and for actions by homeomorphisms of topological spaces
(definition just after Lemma \ref{defsf}), some comments are in order.
Indeed, if a group $\Gamma$ acts on a tree $T$ in a minimal way,
its action on $\partial T$ \emph{need not be minimal}.
A first example is that of the standard action of $\mathbf Z$ on a linear tree;
a second example is the action of a Baumslag-Solitar group $BS(1,n)$
on its Bass-Serre tree $T$, which  is minimal;
but  the corresponding action on the infinite boundary $\partial T$
has a fixed point (see  Subsection \ref{ssHNN} below).
However:

\begin{Prop}
\label{minimal}
Let $\Gamma$ be a group 
which acts  on a tree~$T$.
Assume that the action is  strongly hyperbolic  and minimal.
\par
Then the action of $\Gamma$ on  $\partial T$ is minimal.
%
\end{Prop}
 
\noindent 
\emph{Proof.}~Let $T_0$ be the union of the axis $T_{\gamma}$
over all hyperbolic elements $\gamma \in \Gamma$.
Then $T_0$ is a subtree of $T$ 
and $T_0$ is clearly $\Gamma$-invariant.
Thus $T_0 = T$ by minimality of the action on $T$.
Hence the set
\begin{equation*}
L_{\Gamma} \,  \Doteq \, 
\{ \eta \in \partial T \hskip.1cm \vert \hskip.1cm
\eta = \alpha(\gamma) \hskip.2cm \text{for some hyperbolic} \hskip.2cm
\gamma \in \Gamma \} 
\end{equation*}
is dense in $\partial T$.
\par

The conclusion follows by the argument 
of the proof of Claim (i) of Proposition \ref{Powersgeneral}.
\hfill $\square$

\vskip.2cm

\emph{Remark.}
Proposition \ref{minimal} 
shows that, under appropriate hypothesis,
minimality on $T$ implies minimality on $\partial T$.
Here is a kind of converse:
\par

\emph{
Let $T$ be an infinite tree without vertices of degree one.
If a group $\Gamma$ acts on $T$ in such a way that
its action on $\partial T$ is minimal, 
and if $\Gamma$ contains at least one hyperbolic element,
then the action on $T$ is also minimal.}
\par
Indeed, suppose that the action of $\Gamma$ on $T$ is not minimal,
so that $T_0$ defined as above is a non-empty proper subtree of $T$.
Since $T$ does not have any vertex of degree $1$,
there exists a ray in $T$, disjoint from $T_0$, 
which defines a point $\xi \in \partial T$;
denote by $\left( x_n \right)_{n \ge 0}$ the vertices of this ray,
and by $e$ the first edge of this ray,
with $s(e) = x_0$ and $t(e) = x_1$.
Then the shadow $(\partial T)_e$ is a neighbourhood of $\xi$ in $\partial T$,
disjoint from the $\Gamma$-invariant closed subset $\partial T_0$.
Hence the action of $\Gamma$ on $\partial T$ is not minimal.

\vskip.2cm

The next proposition is a restatement of  well-known facts from the literature.
See for example \cite[Proposition 2]{PaVa--91}
and \cite[Proposition 7.2]{Bass--93}.

\begin{Prop}
\label{transverseh}
Let $T$ be an infinite tree which is not a linear tree.
Let $\Gamma$ be a group which acts on $T$, minimally 
and such that there is no point in $\partial T$ fixed by $\Gamma$.
\par

Then the action of $\Gamma$ on $T$ is strongly hyperbolic.
\end{Prop}

\noindent 
\emph{Proof.}~By minimality, $\Gamma \ne \{1\}$ fixes
no vertex and no pair of adjacent vertices.
By Proposition 3.4 of \cite{Tits--70},
this implies that $\Gamma$ contains at least
one hyperbolic element, say $\gamma_1$,
and that $\partial T \ne \emptyset$.
Denote by $\alpha_1 \in \partial T$ the source of $\gamma_1$
and by $\omega_1$ its sink.

By minimality and
by the proof of Corollary 3.5 of \cite{Tits--70},
$T$ is the union of the axis 
of the hyperbolic elements of $\Gamma$.
As $T$ is infinite and is not a linear tree,
it follows that $\partial T$ is infinite.
The conclusion follows by the  argument 
of the proof of Claim (iii) of Proposition \ref{Powersgeneral}.
\hfill $\square$

\vskip.2cm

When specialised to trees, Proposition \ref{Powersgeneral}
has the following consequence. 

\begin{Cor}
\label{Powerspart}
Let $\Gamma$ be a group which acts on a countable tree $T$. 
Assume that the action is
slender, strongly hyperbolic,  and minimal.
\par
Then $\Gamma$
is a strongly Powers group.
\end{Cor}

\noindent \emph{Proof.}
The action of $\Gamma$ on $\partial T$ is strongly faithful 
by Corollary \ref{slender}, 
and minimal by Proposition \ref{minimal}; 
moreover, we have $\overline{L_{\Gamma}} = \partial T$.
Thus the hypothesis of Claims (i) and (ii) of
Proposition \ref{Powersgeneral} 
are satisfied. 
\par

Let $N \ne 1$ be a normal subgroup of $\Gamma$.
The subgroup $N$ has no fixed vertex in $T$.
Otherwise, since the set $V(T)^N$ is $\Gamma$-invariant,
it would coincide with $V(T)$ by $\Gamma$-minimality, 
and $\Gamma$ could not be faithful (a fortiori not slender) on $T$.
The same argument shows that $N$ has 
no fixed pair of adjacent vertices in $T$
and (because $\Gamma$ is faithful and minimal on $\partial T$)
no fixed boundary point in $\partial T$.
It follows that $N$ contains hyperbolic elements,
by Proposition 3.4 of \cite{Tits--70}.
Thus the hypothesis of Claim (iii) of
Proposition \ref{Powersgeneral} 
are also satisfied. 
\hfill $\square$

\section{\textbf{On the action of the fundamental group 
of a graph of groups on the corresponding Bass-Serre tree
}}
\label{sectionBS}

Recall that a \emph{graph of groups} 
$\mathbf G = (G,Y)$ consists of
\begin{itemize}
\item[$\circledcirc$]
a non-empty connected graph $Y$,
\item[$\circledcirc$]
two families of groups $\left( G_y \right)_{y \in V(Y)}$
and $\left( G_e \right)_{e \in E(Y)}$,
with $G_{\overline{e}} = G_e$ for all $e \in E(Y)$,
\item[$\circledcirc$]
a family of monomorphisms $\varphi_e : G_e \longrightarrow G_{t(e)}$,
for $e \in E(Y)$.
\end{itemize}
An \emph{orientation} of $Y$ is a subset  $E^+(Y)$ of $E(Y)$
containing exactly one of $e,\overline{e}$ for each $e \in E(Y)$;
we denote by $E^-(Y)$ the complement of $E^+(Y)$ in $E(Y)$.
\par

A graph of groups $\mathbf G = (G,Y)$ gives rise to
the \emph{fundamental group} $\Gamma = \pi_1(G,Y,M)$ of $\mathbf G$
and the \emph{universal cover} $T =  T (G, Y, M, E^+(Y))$,
also called the \emph{Bass-Serre tree of $\mathbf G$},
where $M$ is a maximal tree in $Y$,
and  $E^+(Y)$ an orientation of $Y$;
abusively, $T$ is also called the
\emph{Bass-Serre tree of $\Gamma$}.
Let us recall as follows  part of the standard theory
($\S$~I.5 in \cite{Serr--77}, 
and \cite{Bass--93}).

\begin{itemize}

\vskip.2cm

\item[$\circledcirc$]
(BS-1)
$\Gamma$ has a presentation with
generators the groups $G_y$,  $y \in V(Y)$,
and elements $\tau_e$,  $e \in E(Y)$, 
and with relations
\begin{equation*}
\aligned
\tau_{\overline{e}} \, = \, (\tau_e)^{-1} \hskip1cm 
& \text{for all} \hskip.2cm e \in E(Y) ,
\\
\tau_e^{-1} \varphi_{\overline{e}}(h) \tau_e \, = \, \varphi_{e}(h) \hskip1cm 
& \text{for all} \hskip.2cm e \in E(Y) 
    \hskip.2cm \text{and} \hskip.2cm h \in G_e ,
\\
\tau_e \, = \, 1 \hskip1cm
& \text{for all} \hskip.2cm e \in E(M).
\endaligned
\end{equation*}
\par

Moreover, 
the natural homomorphisms 
\begin{equation*}
G_y \longrightarrow \Gamma
\hskip.5cm \text{and} \hskip.5cm
\mathbf Z \longrightarrow \Gamma , \hskip.2cm k \longmapsto \tau_e^k
\end{equation*}
\emph{are injective} for all $y \in V(Y)$ 
and for all $e \in E(Y)$ with  $e \notin E(M)$.

\vskip.2cm

\item[$\circledcirc$]
(BS-2)
$T$ is a graph with 
\begin{equation*}
V(T) \, = \,  \bigsqcup_{y \in V(Y)} \Gamma / G_y
\hskip.5cm \text{and} \hskip.5cm
E(T) \, = \,  \bigsqcup_{e \in E(Y)} \Gamma / \varphi_e(G_e) .
\end{equation*}

The source map, the terminus map, and the inversion map,
are given by
\begin{equation*}
\aligned
s(\gamma \varphi_e(G_e)) \, &= \, \left\{
\aligned
\gamma G_{s(e)} \hskip.2cm
& \text{if} \hskip.2cm e \in E^+(Y) 
\\
\gamma {\tau_e}^{-1} G_{s(e)} \hskip.2cm 
& \text{if} \hskip.2cm e \notin E^+(Y) ,
\endaligned
\right.
\\
t(\gamma \varphi_e(G_e)) \, &= \, \left\{
\aligned
\gamma \tau_e G_{t(e)} \hskip.7cm
& \text{if} \hskip.2cm e \in E^+(Y) 
\\
\gamma G_{t(e)} \hskip.7cm 
& \text{if} \hskip.2cm e \notin E^+(Y) ,
\endaligned \right.
\\
\overline{\gamma \varphi_e(G_e)} \, &= \,  
\gamma \varphi_{\overline e}(G_{\overline{e}}) ,
\endaligned 
\end{equation*}
for all $\gamma \in \Gamma$ and $e \in E(Y)$.
The natural action of $\Gamma$ is by automorphisms of graphs, 
and without inversions.
\par

Moreover, $T$ \emph{is a tree.}

\vskip.2cm

\item[$\circledcirc$]
(BS-3)
The natural mappings
\begin{equation*}
\aligned
V(T) \longrightarrow V(Y) ,
\hskip.5 cm 
&  \gamma G_y \longmapsto y
\\
E(T) \longrightarrow E(Y) ,
\hskip.5cm
& \gamma \varphi_e(G_e) \longmapsto  e
\endaligned
\end{equation*}
are the constituants of a morphism of graphs 
$p : T \longrightarrow Y$
which factors as an isomorphism 
$\Gamma \backslash T \approx Y$.
\par

\vskip.2cm

\item[$\circledcirc$]
(BS-4)
The sections
\begin{equation*}
\aligned
V(Y) \longrightarrow V(T) ,
\hskip.5 cm 
&  y \longmapsto \widetilde y \Doteq 1G_y
\\
E(Y) \longrightarrow E(T) ,
\hskip.5cm
& e \longmapsto \widetilde e \Doteq 1 \varphi_e(G_e)
\endaligned
\end{equation*}
are such that the stabilizer of $\widetilde y$ 
in $\Gamma$ is isomorphic to $G_y$
for all $y \in V(Y)$, and similarly the stabilizer of $\widetilde e$ in $\Gamma$
is isomorphic to $G_e$ for all $e \in E(Y)$.

The first section $V(M) = V(Y) \longrightarrow V(T)$
and the restriction $E(M) \longrightarrow E(T)$ of the second section
are the constituents of an isomorphism of graph
from $M$ onto a subtree of $T$.

\end{itemize}
Moreover, up to isomorphisms, $\Gamma$ and $T$ 
do not depend on the choices of $M$ and $E^+(Y)$. 
\par
For what we need below, it is important to observe that 
\begin{center}
\emph{the action of $\Gamma$ on $T$ 
need  be neither faithful nor minimal.}
\end{center}
We will often write $G_e$ instead of $\varphi_e(G_e)$;
this is abusive since, though $G_{\overline e} = G_e$, the cosets
$\Gamma / G_e$ and $\Gamma / G_{\overline e}$
are {\it different sets} of oriented edges, 
see the definitions of $E(T)$ and of the change of orientations in (BS-2).

\subsection{Reduction of Bass-Serre trees and minimality
}
\label{reduction}

Let $\mathbf G = (G,Y)$ be a graph of groups.
An edge $e \in E(Y)$ is \emph{trivial}\footnote{Terminological adjustment.
``Trivial'' is  here as in \cite[Page 193]{ScWa--79},
so that $e$ is trivial exactly when at least one of $e, \overline{e}$
is ``directed'' as defined in \cite[Page 1096]{Bass--76}.
``Reduced'' here means the same as ``minimal''  in \cite{ScWa--79}.
In \cite[Page 42]{Bass--93}, a vertex $y \in V(Y)$ is 
\emph{terminal} 
if there is a unique edge $e \in E(Y)$ with $s(e) = y$,
and if $\varphi_{\overline{e}}$
is an isomorphism onto; observe then that $e$ is trivial.
\par
If $Y$ is a segment of length $2$ with two end vertices $x,z$,
a middle vertex $y$, the vertex group $G_y$ isomorphic to the two
edge groups, and $G_x, G_z$ large enough,
then $(G,Y)$ is not reduced but does not have any terminal vertex.
The contraction process $\mathbf G \rightsquigarrow \mathbf G / e$
is a particular case of the process described with more details
in \cite[Section 2]{Bass--76}.
}  
if $s(e) \ne t(e)$ 
and if at least one of 
$\varphi_{\overline{e}}$,
$\varphi_{e}$
is an isomorphism onto.
If $\mathbf G = (G,Y)$ has a trivial edge $e$, 
we can define a new graph $Y/e$
obtained from $Y$ by collapsing $\{e,\overline{e}\}$ to a vertex,
and we can define naturally a new graph of groups
$\mathbf G /e$,
with fundamental group isomorphic to that of $\mathbf G$.
Say that $\mathbf G$ is 
\begin{itemize}
\item[$\bullet$]
\emph{reduced}
if it does not contain any trivial edge.
\end{itemize}
\par

Let  $X$ be a connected subgraph of $Y$.
The corresponding \emph{subgraph of groups} $\mathbf F = (F,X)$
is defined by 
$F_x = G_x$ for all $x \in V(X)$,
and $F_e = G_e$ for all $e \in E(X)$,
with the inclusion $F_e \longrightarrow F_{t(e)}$
being precisely $\varphi_e$
for all $e \in E(X)$.
\par

The next proposition collects observations from papers by Bass; 
we choose a maximal tree $M$ of $Y$ containing a maximal tree $L$ of $X$,
and an orientation $E^+(Y)$  of $Y$ containing an orientation $E^+(X)$ of $X$.

\begin{Prop}
\label{Bass1.14+2.15+2.4}
Let $\mathbf G = (G,Y)$ be a graph of groups,
and let the notation be as above.
\begin{itemize}
\item[(i)]
If $\mathbf F = (F,X)$ is a subgraph of groups of $\mathbf G$, 
the fundamental group $\pi_1(F,X,L)$ is isomorphic
to a subgroup of $\pi_1(G,Y,M)$
and the  universal cover $T (F,X,L,E^+(X))$
can be identified with a subtree of $T (G,Y,M,E^+(Y))$.
\par
In case $\mathbf G$ is, moreover, finite and reduced,
and if $X$ is a proper subgraph of $Y$,
then $\pi_1(F,X,L)$ is a proper subgroup of $\pi_1(G,Y,M)$.

\item[(ii)]
If $Y$ is finite,
there exist a reduced graph of group $\mathbf H = (H,Z)$
and a \emph{contraction}
$\mathbf G = (G,Y) 
\longrightarrow \mathbf H = (H,Z)$
which induces an isomorphism 
from  the fundamental group of $\mathbf G$ 
onto that of $\mathbf H$.
\end{itemize}
\end{Prop}

\noindent 
\emph{Proof.}~
For (i), see
Items 1.14 and 2.15 in \cite{Bass--93}.
For (ii), see
Proposition 2.4 in \cite{Bass--76}.
\hfill $\square$
\par
[Claim (ii) need not hold when $Y$ is infinite;
see the discussion in Section~7 of \cite{ScWa--79}.]

\vskip.2cm

The next proposition follows from the proof of Corollary 3.5 in \cite{Tits--70},
and is also part of Proposition 7.12 in \cite{Bass--93}.

\begin{Prop}
\label{minimalG}
Let $\mathbf G = (G,Y)$ be a graph of groups,
and let $\Gamma, T$ be as above;
we assume for simplicity that 
the diameter of the underlying graph $Y$ is finite.
\par

The action of $\Gamma$ on $T$ is minimal
if and only if $Y$ is reduced.
\end{Prop}

\subsection{Small Bass-Serre trees
}
\label{sssmall}

For the next proposition,
compare for example with Theorem 6.1 in \cite{Bass--76}.

\begin{Prop}
\label{finitelinear}
Let $\mathbf G = (G,Y)$ be a graph of groups,
with universal cover $T$;
we assume  that $\mathbf G$ is reduced. Then:
\begin{itemize}

\item[(i)]
$T$ is finite
if and only if $T$ is reduced to one vertex, 
if and only if $Y$ is reduced to one vertex.

\item[(ii)]
$T$ has no vertex of degree $1$.

\item[(iii)]
If $T$ is infinite, 
it does not have any vertex fixed by $\Gamma$.

\item[(iv)]
$T$ is a linear tree
if and only if  $Y$ is 
\subitem{--}
either a segment
with two vertices and one pair of edges $\{e, \overline{e}\}$,
and  with 
$[G_{s(e)} : \varphi_{\overline{e}}(G_e)] = [G_{t(e)} : \varphi_e(G_e)] = 2$
(case of a degenerate non-trivial amalgam, 
see Subsection \ref{ssamalgamated}), 
\subitem{--}
or a loop with one vertex $y$ and 
one pair of edges $\{e, \overline{e}\}$,
and with 
$\varphi_e(G_e) = G_y = \varphi_{\overline{e}}(G_e)$
(case of a semi-direct product $G_y \rtimes_{\theta} \mathbf Z$,
see Subsection \ref{ssHNN}). 

\item[(v)]
$T$ has a pending ray if and only if it is a linear tree.

\end{itemize}
\end{Prop}

\noindent 
\emph{Proof.}~
Let $\mathbf G = (G,Y)$ be a reduced graph of groups.
If $Y$ is not reduced to one vertex,  it contains 
either a segment of length one which is not trivial, or a loop.
In both cases, $T$ contains a linear subgraph $X$,
by Proposition \ref{Bass1.14+2.15+2.4},
and $\Gamma$ contains an element $\gamma$
which leaves $X$ invariant and which induces on $X$
a hyperbolic translation.
Moreover, any edge of $T$ is contained
in a linear subgraph.
Claims (i) to (iii)  follow.

Suppose that $Y$ contains an edge $e$ with $s(e) \ne t(e)$;
suppose moreover  either that at least one of
$[G_{s(e)} : \varphi_{\overline{e}}(G_e)]$, $[G_{t(e)} : \varphi_e(G_e)]$ 
is at least $3$,
or that $Y$ contains at least one other edge than $e$ and $\overline{e}$.
Then $T$ has an abundance of vertices of degrees at least $3$,
and therefore $T$ does not have pending rays.
Similarly,
 if $Y$ contains a vertex $y$ with at least two loops incident to $y$, 
 then $T$ has an abundance of vertices of degrees at least $4$
 and  $T$ does not have pending rays.
This shows Claims (iv) and~(v).
\hfill $\square$

\section{\textbf{The two standard examples
}}
\label{section5}

\subsection{Amalgamated free products
}
\label{ssamalgamated}

In this  case,
the underlying graph $Y$ is a segment  of length one,
with two vertices $x,y$ 
and one pair of edges $\{e,\overline{e}\}$.
The edge group $G_e$ can be identified 
to a subgroup of both $G_x$ and $G_y$.
The fundamental group of $\mathbf G = (G_x, G_y, G_e, Y)$ 
is the free product with amalgamation, 
or, for short, the \emph{amalgam}
$
\Gamma \, = \, G_x \ast_{G_e} G_y .
$
From now on in this subsection, 
we write $A, B, C$  instead of $G_x, G_y, G_e$, 
so that in particular
\begin{equation*}
\Gamma \, = \, A \ast_C B 
\hskip.5cm \text{acts on its Bass-Serre tree} \hskip.2cm T .
\end{equation*}
\par

The edge set of the universal cover $T$ of $\mathbf G$
consists of two copies of $\Gamma / C$
exchanged by the involution
$e \longmapsto \overline{e}$,
say $E(T) = E^+(T) \sqcup E^-(T)$ 
with $E^+(T) = \Gamma / C$ and $E^-(T) = \overline{\Gamma / C}$.
The vertex set of $T$
is the disjoint union $\Gamma / A \sqcup \Gamma / B$.
The source and terminus mappings
are defined to be the canonical projections
\begin{equation*}
\aligned
s : \Gamma / C \longrightarrow \Gamma / A,
\hskip.5cm &\gamma C \longmapsto \gamma A
\\
t : \Gamma / C \longrightarrow \Gamma / B,
\hskip.5cm &\gamma C \longmapsto \gamma B
\endaligned
\end{equation*}
and
\begin{equation*}
\aligned
s : \overline{\Gamma / C} \longrightarrow \Gamma / B,
\hskip.5cm &\overline{\gamma C} \longmapsto \gamma B
\\
t : \overline{\Gamma / C} \longrightarrow \Gamma / A,
\hskip.5cm &\overline{\gamma C} \longmapsto \gamma A .
\endaligned
\end{equation*}
In particular, the tree $T$ is bipartite regular,
with one class of vertices of degree  $[A : C]$
and the other class of degree $[B : C]$.
The action of $\Gamma$ 
has two orbits on the vertex set;
$\Gamma$ acts transitively
on the orientation $E^+(T)$,
or equivalently on the set of geometric edges of $T$.
\par

The \emph{kernel} of such an amalgam is the subgroup
\begin{equation}
\label{Eqkeramalgam}
\ker(A \ast_C B)  \, = \, 
\bigcap_{\gamma \in \Gamma} \gamma^{-1} C \gamma 
\end{equation}
of $C$, namely the largest subgroup of $C$
which is normal in both $A$ and $B$.
An amalgam is
\begin{itemize}
\item[$\bullet$]
\emph{faithful} if its kernel is reduced to $\{1\}$,
equivalently if the action of $\Gamma$ on $T$ is faithful.
\end{itemize}
With the notation of Claim (i) of Theorem~\ref{MainResult}
observe that
\begin{equation*}
\ker(A \ast_C B)  \, = \, 
\bigcap_{\ell \ge 0} C_\ell \, \subset \, C_k \, \subset \, C_0 \, = \, C
\hskip.5cm \text{for all} \hskip.2cm k \ge 0 .
\end{equation*}
The condition $C_k = \{1\}$ has the following geometrical interpretation:
for $\gamma \in \Gamma$, 
if there exists $e \in E(T)$ such that the $k$-neighbourhood
\begin{equation*}
\mathcal V_k(e) \, = \,
\{ x \in V(T) \hskip.1cm \vert
\min \{ d(x, s(e)), \hskip.1cm d(x,t(e)) \} \le k \}
\end{equation*}
is pointwise fixed by $\gamma$, then $\gamma = 1$.
In particular, if $C_k = \{1\}$ for some $k \ge 0$, 
then the action of $\Gamma$ on $T$ is faithful.
\par

The amalgam is 
\begin{itemize}
\item[$\bullet$] 
\emph{non-trivial} if $A \ne C \ne B$ 
(equivalently if $\mathbf G$ is reduced), and
\item[$\bullet$]
\emph{non-degenerate} if moreover
at least one of the indices $[A : C]$, $[B : C]$
is strictly larger than $2$.
\end{itemize}
 From the definition of the universal cover $T$ of $\mathbf G$,
the amalgam is trivial if and only if 
the diameter of $T$ is finite
(and this occurs if and only if the diameter of $T$ is at most $2$).
Also, the amalgam is non-trivial and degenerate
if and only if $T$ is  a linear tree.

\vskip.2cm

\noindent
\emph{Remarks.}
(i)
A non-trivial amalgam which is faithful is 
\emph{a fortiori} non-degenerate,
unless it is the infinite dihedral group.
\par

(ii)
The condition for $\Gamma$ to act faithfully on $T$
\emph{implies} that $\Gamma$ is icc
(as defined in the introduction);
see \cite[Corollary 2]{Corn--09}.
But the converse does not hold. 
\par

(iii)
For example, if $m,n \ge 2$ are two coprime integers,
the torus knot group 
$\langle a,b \hskip.1cm \vert \hskip.1cm a^m = b^n \rangle
=
\langle a \rangle \ast_{\langle a^m = b^n \rangle} \langle b \rangle$
\emph{does not} act faithfully on its Bass-Serre tree;
its kernel is infinite cyclic, generated by $a^m = b^n$.
Any non-trivial free product $A \ast B$ acts faithfully on its Bass-Serre tree.
\par

(iv)
Let $T$ be a regular tree of some degree $d \ge 3$.
Let $\Gamma$ denote the group of all automorphisms $\gamma$ of $T$
which are either elliptic or hyperbolic with an even translation length;
it is a subgroup of index $2$ in $\operatorname{Aut}(T)$,
and simple \cite{Tits--70}.
Choose an edge $e \in V(T)$.
Denote by $A$ [respectively $B$, $C$] the pointwise stabiliser in $\Gamma$
of $s(e)$ [respectively of $t(e)$, of $\{s(e),t(e)\}$],
so that $\Gamma = A \ast_C B$
\cite[Section~I.5.4, Theorem~13]{Serr--77}.
We leave it to the reader to check that
(as already stated in the introduction), 
for any $k \ge 0$,
the group $C_k$ is the pointwise stabiliser of
the neighbourhood $\mathcal V _k(e)$ defined in the introduction,
and that
\begin{equation*}
\aligned
C \, = \, C_0 \,  \supsetneqq \,   C_1 \,  &\supsetneqq \, \cdots \, 
\supsetneqq  \, C_k  \,  \supsetneqq  \, C_{k+1} \,  \supsetneqq  \, \cdots 
\\
&\supsetneqq  \bigcap_{k \ge 0}C_k 
\, = \, \ker(A \ast_C B) \, = \, \{1\} 
\endaligned
\end{equation*} 
is a \emph{strictly} decreasing infinite sequence of subgroups of $C$.
As already noted in Remark \ref{Remarks}.v, 
the action of $\Gamma$ on $T$ is not slender.
\par

Observe that this group $\Gamma$ is not countable,
but that it contains dense countable subgroups
giving rise similarly to strictly decreasing infinite sequences
of $C_k$'s.
We are grateful to Laurent Bartholdi for suggesting these examples.
\par

\emph{ We do not know if these groups are C$^*$-simple.}
\par

(v)
For an integer $k \ge 1$, recall that $\Gamma = A \ast_C B$
is \emph{$k$-acylindrical} if, 
whenever $\gamma \in \Gamma$  fixes pointwise
a segment of length $k$ in its Bass-Serre tree,
then $\gamma = 1$ \cite{Sela--97}.
Thus the condition ``$C_k = \{1\}$ for some $k \ge 0$''
is substancially weaker than the condition
``$\Gamma$ is $\ell$-acylindrical for some $\ell \ge 0$''.
The  condition ``$\ker (A \ast_C B) = \{1\}$''
is the weakest of all.

\begin{Prop}
\label{PropAmal}
Let $\Gamma = A \ast_C B$ be an 
amalgam
acting on its Bass-Serre tree $T$ as above.
\begin{itemize}
\item[(i)]
The amalgam is non-trivial if and only if
there is no vertex in $T$ fixed by $\Gamma$,
if and only if the tree $T$ is infinite;
if these conditions hold, then
the action of $\Gamma$ on $T$ is \emph{minimal}.
\end{itemize}
We assume from now on that the amalgam is non-trivial.
\begin{itemize}
\item[(ii)]
The amalgam is non-degenerate if and only if
$T$ is not a linear tree,
if and only if $\partial T$ is perfect,
if and only if
the action of $\Gamma$ on $T$ is \emph{strongly hyperbolic}.
\end{itemize}
We assume from now on that the amalgam is non-degenerate. 
\begin{itemize}
\item[(iii)]
The action of $\Gamma$ on $\partial T$ is minimal.
\item[(iv)]
The action of $\Gamma$ on $T$ is slender
as soon as $C_k = \{1\}$ for some $k \ge 1$.
\end{itemize}
\end{Prop}

\noindent 
\emph{Remark.}~Together with (ii), recall also the following particular case 
of Theorem 6.1 in \cite{Bass--76} :
Let $\Gamma = A \ast_C B$ be a non-trivial amalgam;
then $\Gamma$ does not contain non-abelian free subgroups
if and only if $\Gamma$ is degenerate
and $C$ does not contain non-abelian free subgoups. 

\vskip.2cm

\noindent 
\emph{Proof.}~For (i), 
see Propositions \ref{minimalG} and  \ref{finitelinear}. 
[Alternatively, a direct argument is straightforward.]

\vskip.1cm

(ii) Observe first that, if the amalgam is non-trivial and degenerate,
then $T$ is a linear tree,
so that $\partial T$ has two points, in particular is not perfect,
and the action of $\Gamma$ on $\partial T$
is not strongly hyperbolic.
\par
Suppose now that the amalgam is non-degenerate.
Choose $q,r \in A$ such that the three cosets $C, qC, rC$ in $A/C$
are pairwise disjoint, 
and $s \in B$ with $s \notin C$
(so that $C \cap sC = \emptyset = C \cap s^{-1}C$ in $B/C$).
\par
In $T$, there is a first segment of length $5$, 
of which the $6$ vertices and $5$ oriented edges are,
in ``the'' natural order,
\begin{equation*}
qs^{-1}A, \hskip.2cm  qs^{-1}C, \hskip.2cm
qB, \hskip.2cm qC, \hskip.2cm 
\underline{A}, \hskip.2cm \underline{C}, \hskip.2cm 
\underline{B}, \hskip.2cm \underline{sC}, \hskip.2cm 
\underline{sA}, \hskip.2cm sqC , \hskip.2cm 
sqB,
\end{equation*}
and similarly a second segment of length $5$
with vertices and edges
\begin{equation*}
rs^{-1}A, \hskip.2cm  rs^{-1}C, \hskip.2cm
rB, \hskip.2cm rC, \hskip.2cm 
\underline{A}, \hskip.2cm \underline{C}, \hskip.2cm 
\underline{B}, \hskip.2cm \underline{sC}, \hskip.2cm 
\underline{sA}, \hskip.2cm srC , \hskip.2cm 
srB.
\end{equation*}
These two segments of length $5$
have a common subsegment of length $2$
with vertices and edges underlined above. 
\par

The element $sqsq^{-1} \in \Gamma$ maps 
the first edge of the first segment onto its last edge,
and $srsr^{-1} \in \Gamma$ maps
the first edge of the second segment onto its last edge.
It follows that $sqsq^{-1}$ and $srsr^{-1}$
are hyperbolic elements, with axis sharing exactly 
two geometric edges
(see if necessary Proposition 25 in Subsection I.6.4 of \cite{Serr--77}).
Thus these two hyperbolic elements are transverse.
\par

If the amalgam is non-degenerate,
at least every other vertex in $T$ is of degree at least $3$.
In particular, $T$ does not have any pending ray
(Proposition \ref{finitelinear}), so that $\partial T$ is perfect.

\vskip.1cm

(iii)
If the amalgam is non-degenerate,
the action of $\Gamma$ on $T$ is minimal, by (i).
Hence the action is also minimal on $\partial T$, 
by Proposition \ref{minimal}. 
\par

\vskip.1cm

(iv)
Let $\gamma \in \Gamma$ be such that 
the fixed point set $(\partial T)^{\gamma}$ of $\gamma$ 
on the boundary has non-empty interior.
To finish the proof, 
it suffices to show that, 
if the amalgam is non-degenerate and if $C_k = \{1\}$,
then $\gamma = 1$.
\par

Since $\partial T$ is perfect, $\gamma$ cannot be hyperbolic.
We can therefore assume that $\gamma$ has a fixed vertex $x_0 \in V(T)$.
Let $(x_n)_{n \ge 0}$ be the vertices of a ray
starting from $x_0$ and representing a boundary point $\xi$
in the interior of $(\partial T)^{\gamma}$.
By Remark  \ref{Remarks}.iv, 
the vertex $x_n$ is fixed by $\gamma$ for all $n \ge 0$.
Since $\xi$ is in the interior of $(\partial T)^{\gamma}$,
there is an edge  in the ray  $(x_n)_{n \ge 0}$,
say $d$ from $x_m$ to $x_{m+1}$,
such that $(\partial T)_d \subset (\partial T)^{\gamma}$.
\par

Let $U$ be the subtree of $T$ of which the vertices
belong to rays with first two vertices $x_m, x_{m+1}$.
Since rays in this tree represent 
boundary points in $(\partial T)^{\gamma}$,
the same remark as above implies that
$\gamma$ fixes all vertices and all edges in $U$.
Choose an edge $e \in E(U)$ such that 
all vertices at distance at most $k$ from $s(e)$ or $t(e)$
are in $V(U)$; these vertices are  fixed by $\gamma$.
Choose moreover $\delta \in \Gamma$ such that $e = \delta C$;
we can assume that $s(e) \in \Gamma / A$ and $t(e) \in \Gamma / B$.
\par

Let us first assume for simplicity that $k = 1$.
Choose  transversals $R \subset A$ and $S \subset B$ such that 
$A = \bigsqcup_{r \in R} rC$ and $B = \bigsqcup_{s \in S} sC$
(disjoint unions).
Since $\gamma$ fixes the edges $\delta rC$ and $\delta sC$,
we have $\gamma \in \delta rC(\delta r)^{-1}$ 
and $\gamma \in \delta sC(\delta s)^{-1}$ 
for all $r \in R$ and $s \in S$,
namely
\begin{equation*}
\aligned
\gamma \, &\in \,
\big( \bigcap_{r \in R} \delta rCr^{-1}\delta^{-1} \big)  
\cap
\big( \bigcap_{s \in R} \delta sCs^{-1}\delta ^{-1} \big)
\\
\, &= \,
\delta \left( 
\big( \bigcap_{a \in A} aCa^{-1}  \big)
\cap
\big( \bigcap_{b \in B} bCb^{-1} \big)
\right) \delta ^{-1} 
\\
\, &= \,
\delta C_1 \delta^{-1}
\, = \,
\{1\} ,
\endaligned
\end{equation*}
hence $\gamma = 1$.
\par

The argument in the general case, $k \ge 1$, is similar,
and is left to the reader. Hence, in all cases, $\gamma = 1$.
\hfill $\square$

\vskip.2cm

\noindent
\textbf{Proof of Claim (i) of Theorem~\ref{MainResult}.}
We have now to assume that
$A$, $B$, and therefore also 
$C$, $\Gamma$, and $T$, are countable,
because our proof of Proposition \ref{Baire} 
assumes countability. 
Moreover, we assume as in Claim (i)
that the amalgam is non-degenerate 
and that $C_k = \{1\}$.

\par
Proposition \ref{PropAmal} 
shows that the hypothesis of Corollary \ref{Powerspart} 
are satisfied.
Hence the proof of Claim (i) is complete.
\hfill $\square$

\subsection{HNN extensions
}
\label{ssHNN}

In this  case,
the underlying graph $Y$ is a loop,
with one vertex $y$, 
and one pair of edges $\{e,\overline{e}\}$.
The edge group $G_e$ can be identified 
(via $\varphi_{\overline{e}}$) to a subgroup of $G_y$,
and we have a monomorphism $\varphi_e : G_e \longrightarrow G_y$.
From now on in this subsection, 
we write $G, H, \theta$ instead of $G_y, G_e, \varphi_{e}$.
The fundamental group of $\mathbf{G} = (G, H, \theta(H), Y)$ is a
\emph{HNN-extension}, which has the presentation
\begin{equation*}
\Gamma \, = \, HNN(G, H, \theta) \, = \,
\langle G, \tau \hskip.1cm \vert \hskip.1cm
\tau^{-1}h\tau = \theta(h) \hskip.2cm
\forall h \in H \rangle 
\end{equation*}
and which acts on its Bass-Serre tree $T$.
\par

The edge set of  $T$ 
consists of two copies of $\Gamma / H$
exchanged by the involution $e \longmapsto \overline{e}$,
say $E(T) = E^+(T) \sqcup E^-(T)$ 
with $E^+(T) = \Gamma / H$ and $E^-(T) = \overline{\Gamma / H}$.
The vertex set of $T$
is  $\Gamma / G$.
The source and terminus mappings
are 
\begin{equation*}
\aligned
s : \Gamma / H \longrightarrow \Gamma / G , \hskip.5cm
&\gamma H \longmapsto \gamma G
\\
t : \Gamma / H \longrightarrow \Gamma / G , \hskip.5cm
&\gamma H \longmapsto \gamma \tau G .
\endaligned
\end{equation*}
and
\begin{equation*}
\aligned
s : \overline{\Gamma / H} \longrightarrow \Gamma / G , \hskip.5cm
&\overline{\gamma H} \longmapsto \gamma \tau G
\\
t : \overline{\Gamma / H} \longrightarrow \Gamma / G , \hskip.5cm
&\overline{\gamma H} \longmapsto \gamma  G .
\endaligned
\end{equation*}
In particular, the tree $T$ is regular, 
of degree $[G:H] + [G : \theta(H)]$.
The action of $\Gamma$ is transitive  
on the vertex set,
as well as on the orientation $E^+(T)$, 
or equivalently on the set of geometric edges of $T$.
\par

The \emph{kernel} of an HNN-extension is the subgroup
\begin{equation}
\label{EqkerHNN}
\ker(HNN(G,H,\theta))  \, \Doteq \, 
\bigcap_{\gamma \in \Gamma} \gamma^{-1} H \gamma 
\end{equation}
of $H \cap \theta(H)$, namely the largest subgroup of $H \cap \theta(H)$
which is both normal in $G$ and invariant by $\theta$.
An HNN extension is
\begin{itemize}
\item[$\bullet$]
\emph{faithful} if its kernel is reduced to $\{1\}$,
equivalently if the action of $\Gamma$ on $T$ is faithful.
\end{itemize}
With the notation of  Claim (ii) in Theorem~\ref{MainResult},
observe that
\begin{equation*}
\ker(HNN(G,H,\theta)) \, = \,
\bigcap_{\ell \ge 0} H_{\ell} \, \subset \, H_k  \, \subset \, H'_1 \, = \, H \cap \theta(H)
\hskip.5cm \text{for all} \hskip.2cm k \ge 1.
\end{equation*}
The condition $H_k = \{1\}$ has the following geometrical interpretation:
for $\gamma \in \Gamma$, 
if there exists $e \in E(T)$ such that
the $k$-neighbourhood $\mathcal V_k(e)$
is pointwise fixed by $\gamma$, then $\gamma = 1$.
\par

\par

The HNN-extension is 
\begin{itemize}
\item[$\bullet$] 
\emph{ascending} 
if at least one of $H, \theta(H)$ is the whole of $G$,
\item[$\bullet$]
\emph{strictly ascending} 
if exactly one of $H, \theta(H)$ is the whole of $G$, and
\item[$\bullet$]
\emph{non-degenerate} if at least one of $H, \theta(H)$ 
is a proper subgroup of $G$.
\end{itemize}
The HNN-extension is degenerate 
if and only if $T$ is a linear tree,
in which case $\theta$ is an automorphism of $G$
and $\Gamma$ is the corresponding semi-direct product
$G \rtimes_{\theta} \mathbf Z$.

\vskip.2cm

\noindent
\emph{Remarks.}
(i)
An HNN-extension with $H \ne \{1\}$
which is faithful is \emph{a fortiori} non-degenerate.
\par

(ii)
The condition for $\Gamma$ to act faithfully on $T$
\emph{implies} that $\Gamma$ is icc
(but the converse does not hold); 
see  \cite[Example 2.9]{Stal--06}
and \cite[Corollary 4]{Corn--09}.
\par

(iii)
For example, if $m,n$ are  integers, 
the Baumslag-Solitar group
\begin{equation*}
BS(m,n) \, = \, 
\langle \tau,b \hskip.1cm \vert \hskip.1cm \tau^{-1}b^m\tau = b^n \rangle
\, = \,
HNN(b^{\mathbf Z}, b^{m\mathbf Z}, b^{mk} \longmapsto b^{nk})
\end{equation*}
acts faithfully on the corresponding tree
if and only if $\vert m \vert \ne \vert n \vert$.
Indeed, on the one hand,  if $n = \pm m$, then
$H = \langle b^m \rangle$ is clearly 
a normal subgroup of  $\Gamma$;
and, on the other hand,
it is a result of Moldavanskii
that $BS(m,n)$ has an  infinite cyclic normal subgroup
if and only if $\vert m \vert = \vert n \vert$;
see \cite{Mold--91}, 
or the exposition in the Appendix of \cite{Souc--01}.

\begin{Prop}
\label{PropHNN}
Let $\Gamma = HNN(G,H,\theta)$ be an 
HNN-extension
acting on its Bass-Serre tree $T$ as above. Then:
\begin{itemize}
\item[(i)]
There is no vertex in $T$ fixed by $\Gamma$,
the tree $T$ is infinite and the action of $\Gamma$ on $T$ is minimal.
\item[(ii')]
The HNN-extension is non-degenerate if and only if 
$T$ is not a linear tree,
if and only if the space $\partial T$ is perfect.
\end{itemize}
We assume from now on that the extension is non-degenerate.
\begin{itemize}
\item[(ii'')]
The HNN-extension is non-ascending if and only if 
the space $\partial T$ is  without $\Gamma$-fixed point,
if and only if the action of $\Gamma$ on $T$ is strongly hyperbolic.
\end{itemize}
We assume from now on that the extension is non-ascending.
\begin{itemize}
\item[(iii)]
The action of $\Gamma$ on $\partial T$ is minimal.
\item[(iv)]
The action of $\Gamma$ on $T$ is slender
as soon as $H_k = \{1\}$ for some $k \ge 1$.
\end{itemize}
\end{Prop}

\noindent 
\emph{Remark.}~Together with (ii'), recall also the following particular case 
of Theorem 6.1 in \cite{Bass--76} :
Let $\Gamma = HNN(G,H,\theta)$ be a  HNN-extension;
then $\Gamma$ does not contain non-abelian free subgroups
if and only if $\Gamma$ is ascending
and $G$ does not contain non-abelian free subgroups. 

\vskip.2cm

\noindent 
\emph{Proof.}~
For (i) and (iii), see the proof of Proposition \ref{PropAmal}. 

\par

(ii') If the extension is degenerate, $T$ is a linear tree,
so that $\partial T$ has exactly two points.
Otherwise, $T$ is regular of degree has least $3$,
so that $\partial T$ is perfect.

\par

(ii'')
Suppose first that $H = G \supsetneqq \theta(H)$. The mapping 
$E^+(T) \longrightarrow V(T), \hskip.2cm e \longmapsto s(e)$
is a bijection. 
Choose a geodesic in $T$ with vertex set 
$\left( x_p \right)_{p \in \mathbf Z}$
such that $d(x_p, x_{p'}) = \vert p - p' \vert$
for all $p,p' \in \mathbf Z$,
and such that the edge from $x_p$ to $x_{p+1}$ 
lies in $E^+(T)$ for all $p \in \mathbf Z$;
consider the limit $\xi = \lim_{p \to \infty} x_p \in \partial T$.
Consider also a geodesic ray with vertex set
 $\left( y_q \right)_{q \in \mathbf N}$
such that $d(y_q, y_{q'}) = \vert q - q' \vert$
for all $q,q' \in \mathbf N$,
and such that the edge from $y_q$ to $y_{q+1}$ 
lies in $E^+(T)$ for all $q \in \mathbf N$;
set $\eta = \lim_{q \to \infty} y_q \in \partial T$.
We claim that $\eta = \xi$,
from which it follows that $\xi$ is fixed by $\Gamma$.
\par

To prove the claim, it suffices to check that
there exists $b \in \mathbf Z$ and $q \in \mathbf N$
such that $y_q = x_{q+b}$;
indeed, since $E^+(T) \longrightarrow V(T)$, $e \longmapsto s(e)$,
is a bijection, this implies that $y_{q'} = x_{q' + b}$ for all $q' \ge q$,
and therefore $\eta = \xi$.
\par

If one could not find $b,q$ such that $y_q = x_{q+b}$,
there would exist a segment with vertex set
$\left( z_r \right)_{0 \le r \le N}$, with $N \ge 1$,
connecting the geodesic to the ray,
namely with $z_0 = x_{p_0}$ for some $p_0 \in \mathbf Z$
and $z_N = y_{q_0}$ for some $q_0 \in \mathbf N$,
moreover with
$z_r \notin \left\{ x_p \right\}_{p \in \mathbf Z}$ for $r \ge 1$
and 
$z_r \notin \left\{ y_q \right\}_{q \in \mathbf N}$ for $r \le N-1$.
On the one hand, the edge with source $z_{r+1}$ and terminus $z_r$
would be in $E^+(T)$ for all $r$
(as one checks inductively for $r = 0, 1, \hdots$),
on the other hand, the edge with source $z_{r-1}$ and terminus $z_r$
would also be in $E^+(T)$ for all $r$
(as one checks inductively for $r = N, N-1, \hdots$).
But this is impossible, since $E^+(T)$ cannot contain both
some edge $e$ with source $z_{r+1}$ and terminus $z_r$
and the edge $\overline{e}$.
\par

Suppose next that $H \subsetneqq G = \theta(H)$,
so that the mapping 
$E^+(T) \longrightarrow V(T), \hskip.2cm e \longmapsto t(e)$
is a bijection.
An analogous argument shows that there exists a point
in $\partial T$ fixed by $\Gamma$.
\par

It follows also that, if the HNN-extension is ascending,
there cannot exist two transverse hyperbolic elements in $\Gamma$.
\par

Suppose now that the extension is non-ascending.
Choose $r,s \in G$ with $r \notin H$ and $s \notin \theta(H)$.
\par
In $T$, there are two segments of length $2$,
with vertices and edges respectively
\begin{equation*}
\begin{matrix}
\tau^{-1}G & \tau^{-1}H & G & rH & r\tau G ,
\\
s^{-1}\tau^{-1}G & s^{-1}\tau^{-1}H & G & H & \tau G ,
\end{matrix}
\end{equation*}
sharing just one vertex, $G$.
It follows that $r\tau$ and $\tau s$ are two elements of $\Gamma$
which are hyperbolic and transverse,
with axis having in common the unique vertex $G$.

\par

(iv) The proof of this claim is a somewhat tedious variation on that
of Proposition \ref{PropAmal} 
and is left to the reader.
\hfill $\square$

\vskip.2cm

\noindent
\textbf{Proof of Claim (ii) of Theorem~\ref{MainResult}.}
We have now to assume that
$G$, and therefore also $H$, $\Gamma$, and $T$, are countable,
because our proof of Proposition \ref{Baire} assumes countability. 
Moreover, we assume as in Claim (ii)
that the HNN-extension is non-ascending and that $H_k = \{1\}$.
\par
Proposition \ref{PropHNN} 
shows that the hypothesis of Corollary \ref{Powerspart} 
are satisfied.
Hence the proof of Claim (ii) is complete.
\hfill $\square$

\begin{Lemma}
\label{LemmaBS}
Consider the Baumslag-Solitar group
\begin{equation*}
\Gamma = \langle \tau, b \hskip.1cm \vert \hskip.1cm 
\tau^{-1}b^m\tau = b^n \rangle
\end{equation*}
acting on its Bass-Serre tree $T$. 
Denote by $H$ the cyclic subgroup of $\Gamma$ generated by $b^m$.
Set
\begin{equation*}
K_+ \, = \, \bigcap_{p \ge 0} \tau^{-p} H \tau^{p} 
\hskip.5cm \text{and} \hskip.5cm
K_- \, = \, \bigcap_{p \ge 0} \tau^{p} H \tau^{-p} .
\end{equation*}

(i) We have:
\begin{itemize}
\item[$\circ$] 
if $\vert m \vert = \vert n \vert$, then $K_+ = K_- = H$,
\item[$\circ$] 
if $n = \pm am$ for some $a \ge 2$, then $K_+ = \{1\}$ and $K_- = H$,
\item[$\circ$] 
if $m = \pm an$ for some $a \ge 2$, then $K_+ = H$ and $K_- = \{1\}$,
\item[$\circ$] 
in all other cases,  $K_+ = K_- = \{1\}$.
\end{itemize}
\par

(ii) Assume that $\vert m \vert \ne \vert n \vert$.
For any $k \in \mathbf Z$, $k \ne 0$, the automorphism $g = b^k$
of $T$ is elliptic and slender.
\par

(iii) If $\vert m \vert \ne \vert n \vert$, the action of $\Gamma$ on $T$ is slender.
\end{Lemma}

\noindent 
\emph{Proof.}~(i) 
Denote by $G$ the cyclic subgroup of $\Gamma$ generated by $b$,
so that $H \subset G$ and $[G:H] = m$.
The normal form theorem for HNN extensions
(Theorem 2.1 of Chapter IV in \cite{LySc--77})
implies that
\begin{equation*}
\aligned
\tau^{-1} b^\ell \tau \, &\notin \, G
\hskip.5cm \text{if} \hskip.2cm \ell \notin m \mathbf Z,
\\
\tau^{-1} b^{km} \tau \, &= \, b^{kn}
\hskip.5cm \text{for all} \hskip.2cm k \in \mathbf Z,
\\
\tau b^\ell \tau^{-1} \, &\notin \, G
\hskip.5cm \text{if} \hskip.2cm \ell \notin n \mathbf Z,
\\
\tau b^{kn} \tau^{-1} \, &= \, b^{km}
\hskip.5cm \text{for all} \hskip.2cm k \in \mathbf Z .
\endaligned
\end{equation*}
Claim (i) follows.

\vskip.1cm

For (ii) and (iii), we assume that $K_- = \{1\}$.
The case with $K_+ = \{1\}$ will follow,
since $HNN(G,H,\theta) \approx HNN(G, \theta(H), \theta^{-1})$.
Let  $x_0$ denote the vertex $G = 1G \in \Gamma / G = V(T)$.
Observe that $x_0$ is fixed by $g$;
indeed, the isotropy subgroup 
$\{\gamma \in \Gamma \hskip.1cm \vert \hskip.1cm \gamma x_0 = x_0 \}$
coincides with $G$.

\vskip.1cm

(ii) 
It suffices to show the following claim:
$(\partial T)_e \not\subset (\partial T)^g$ 
for any $e \in E(T)$ with $d(x_0,s(e)) < d(x_0,t(e))$.
If $g(t(e)) \ne t(e)$, then 
$g\left( (\partial T)_e \right) \cap (\partial T)_e = \emptyset$,
and the claim is obvious; we can therefore assume that $g(t(e)) = t(e)$.
We distinguish two cases, 
depending on the $\Gamma$-orbit of~$e$.
\par

Suppose first that $e \in \Gamma / H$, namely that
there exists $\gamma \in \Gamma$ such that
$s(e) = \gamma \tau^{-1} (x_0)$ and $t(e) = \gamma (x_0)$.
Then $\left( \gamma \tau^p G \right)_{p \in \mathbf N}$
are the vertices of a ray $\rho$ in $T$ starting at $t(e)$,
with $d(\gamma \tau^{p} G, \gamma \tau^{q}G) = \vert p-q \vert$
for all $p,q \in \mathbf N$,
and extending the segment from $x_0$ to $t(e)$. 
For $p \ge 1$, the vertex $\gamma \tau^p G$ is fixed by $g$
if and only if
$g  \in  \gamma G \gamma^{-1} \cap \gamma \tau^p G \tau^{-p} \gamma^{-1}$,
namely if and only if
$\gamma^{-1} g \gamma  \in  G \cap  \tau^p G \tau^{-p}$.
It follows from (i) that there exists an edge $f$ in the ray $\rho$,
with $d(x_0, s(f)) < d(x_0,t(f))$, 
such that $g(s(f)) = s(f)$ and $g(t(f)) \ne t(f)$,
and consequently such that 
$g\left( (\partial T)_f \right) \cap (\partial T)_f = \emptyset$.
Since $(\partial T)_f  \subset (\partial T)_e$,
this implies $(\partial T)_e \not\subset (\partial T)^g$.
\par

Suppose now that $e \in \overline{\Gamma / H}$, namely that
there exists $\gamma \in \Gamma$ such that
$s(e) = \gamma \tau (x_0)$ and $t(e) = \gamma (x_0)$.
Observe that $H \ne G$, otherwise $K_- = G$ could not be $\{1\}$.
Choose $u \in G$ with $u \notin H$,
so that $\overline{\gamma u H}$ is an edge with source
$\gamma u \tau G \ne s(e) = \gamma \tau G$
and with terminus $\gamma G = t(e)$. 
Then $(\gamma u \tau^p G)_{p \in \mathbf N}$
are the vertices of a ray in $T$ 
starting at $\gamma G = \gamma u G = t(e)$,
and the argument of the previous case carries over.

\vskip.1cm

(iii) Let $\gamma \in \Gamma$, $\gamma \ne 1$,
be an elliptic automorphism of $T$.
Choose $\delta \in \Gamma$ such that
the vertex $\delta G$ is fixed by $\gamma$.
Then $g  \Doteq \delta^{-1} \gamma \delta$ fixes $x_0$,
so that $g = b^k$, as in Claim (ii).
We have $(\partial T)^{\gamma} = \delta \left( (\partial T)^g \right)$;
in other words, the subspace $(\partial T)^{\gamma}$ of $\partial T$
is the image by the homeomorphism $\delta$ of $\partial T$
of the subspace  $(\partial T)^g$, without interior points by (ii).
Hence $\gamma$ is slender.
\hfill $\square$

\vskip.2cm 

\noindent
\textbf{Proof of Claim (iii) of Theorem~\ref{MainResult}.}
If $\min \{ \vert m \vert, \vert n \vert \} = 1$,
the Baumslag-Solitar group $BS(m,n)$ is solvable,
in particular amenable, and therefore it is not C$^*$-simple.
If $\vert m \vert = \vert n \vert$, the group $BS(m,n)$
contains an infinite cyclic normal subgroup
(as already noted in Remark (iii) 
just before Proposition \ref{PropHNN}), 
so that $BS(m,n)$ is not C$^*$-simple. 
\par
Let us assume that 
$\min \{ \vert m \vert, \vert n \vert \} \ge 2$
and that $\vert m \vert \ne \vert n \vert$.
Observe that, say in the case of $m$ and $n$ coprime to simplify the discussion,
we have  $G = \mathbf Z$ and
\begin{equation*}
\aligned
H_0 = m \mathbf Z \supsetneqq H_1 = (mn)^2  \mathbf Z \supsetneqq  \cdots
&\supsetneqq  H_k = (mn)^{2k} \mathbf Z \supsetneqq  \cdots 
\\
&\supsetneqq  \bigcap_{k \ge 0}H_k = \ker(BS(m,n)) = \{0\} ,
\endaligned
\end{equation*} 
so that we cannot apply 
Claim (ii) in Theorem~\ref{MainResult}, 
even though the action of $BS(m,n)$ on $T$ is faithful
(by this same Remark (iii)).
\par

However,  the action of $BS(m,n)$ on
its Bass-Serre tree $T$ 
is slender, by Lemma \ref{LemmaBS}, 
strongly hyperbolic and minimal, by Proposition \ref{PropHNN}. 
Hence, the hypothesis of Corollary \ref{Powerspart} 
are satisfied,
and the proof of Claim (iii) of Theorem~\ref{MainResult}
for $BS(m,n)$ is complete.
\par
The proof of the more general case, which is similar, is left to the reader.
\hfill $\square$

\subsection{Remark
}
\label{ssremarks2}

Let $\mathbf G = (G,Y)$ be a  graph of groups.
If $Y$ is a tree, the analysis of the fundamental group of $\mathbf G$
can be done essentially as in Subsection \ref{ssamalgamated}. 
If $Y$ is not a tree, one can first analyze the graph of groups
related to a maximal subtree of $Y$,
and then analyze the contributions of the remaining edges
as in Subsection \ref{ssHNN}.  
\par

On several occasions in the past, 
including in \cite{Harp--85} and \cite{BeHa--86}, 
the first author has been wrong
in dealing with  elliptic automorphisms of a tree $T$,
concerning their fixed point sets in $V(T)$ and in $\partial T$.
As a consequence of this,
Claims (d) and (e) of Theorem 5 in \cite{BeHa--86} 
are not correct as stated, as pointed out by Yves Stalder in \cite{Stal--06}
(many thanks to him).
The notion of \emph{slender} automorphism
provides some way to fix
at least part of the confusion.

\vskip.4cm
\begin{center}
\textbf{Second part: on fundamental groups of $3$-manifold}
\end{center}

Manifolds which appear below are assumed to be \emph{connected},
but for boundaries and a few other \emph{explicit} exceptions. 
If $M$ is a manifold, we denote by $\partial M$ its boundary.
We use $\approx$ to indicate both 
a diffeomorphism of manifolds
and an isomorphism of groups.
\par

We also use rather standard notation for particular manifolds,
with the dimension in superscript:
$\mathbf S^n$ for spheres, 
$\mathbf P^n$ for real projective spaces,
$\mathbf I$, $\mathbf D^2$, and $\mathbf B^3$ 
for the interval, the $2$-disc, and the $3$-ball,
$\mathbf T^2$ and $\mathbf K^2$ 
for the $2$-torus and the Klein bottle.
Recall that there are two disc bundles over the circle:
the trivial bundle, which is the direct product $\mathbf S^1 \times \mathbf D^2$,
is the \emph{solid torus};
the non-trivial bundle, denoted by 
$\mathbf S^1 \widetilde \times \mathbf D^2$,
is a non-orientable $3$-manifold with boundary a Klein bottle,
called the \emph{solid Klein bottle}.
Similarly, there are two $\mathbf I$-bundles over $\mathbf K^2$:
one is trivial and with non-orientable total space,
the other is non-trivial with orientable total space;
the second
is denoted here by $\mathbf K^2 \widetilde{\times} \mathbf I$.

For background on $3$-manifolds, 
we refer to  \cite{Hemp--76}, \cite{JaSh--79}, 
\cite{Scot--83},  \cite{Thur--97}, and \cite{Bona--02}.

\section{\textbf{Some particular cases of Theorem 2.iv
}}
\label{sectionparticularcases}

Before stating the main result of this section,
we review the following standard definitions:
\begin{itemize}
\item[$\bullet$]
Seifert manifolds,
\item[$\bullet$]
Sol-manifolds,
\item[$\bullet$]
hyperbolic manifolds,
\item[$\bullet$]
irreducibility,  incompressibility, and $\partial$-irreducibility,
\item[$\bullet$]
essential tori.
\end{itemize}
Readers comfortable with these notions 
should go to Proposition \ref{examples-3-manif}.

\vskip.2cm

(\ref{sectionparticularcases}.1). 
A \emph{Seifert manifold} is a compact $3$-manifold 
which can be foliated by circles, 
here called \emph{leaves} (but often called \emph{fibers} elsewhere),
such that each leaf has a foliated neighbourhood, 
either a solid torus 
or a solid Klein bottle,
which is finitely covered by the solid torus $\mathbf S^1 \times \mathbf D^2$
with the product foliation $(\mathbf S^1 \times \{z\})_{z \in \mathbf D^2}$,
in such a way that leaves cover leaves.
\par

For example, any compact orientable $3$-manifold 
with finite fundamental group
and without spherical boundary component
is a Seifert manifold.
(This statement, now a result by Perelman's work,
was known as the elliptisation conjecture,
a part of Thurston's geometrisation conjecture.)
Let $M$ be a Seifert manifold 
with infinite fundamental group $\Gamma$
(equivalently: such that $M$ is not covered by a $3$-sphere).
If $M$ is a Seifert manifold, 
a generic leaf generates a normal subgroup of $\Gamma$
which is infinite cyclic \cite[Lemma 3.2]{Scot--83}.

\vskip.2cm

(\ref{sectionparticularcases}.2)
A \emph{Sol-manifold} is a $3$-manifold $M$,
here for simplicity assumed to be without boundary,
such that the universal covering of $M$ 
can be identified with $Sol$
and such that $M$ has a Riemannian structure
for which the universal covering 
$\pi : Sol \longrightarrow M$ is a local isometry.
Here, $Sol$ denotes the $3$-dimensional Lie group
with underlying space $\mathbf R^3$,
with product
\begin{equation*}
(x,y,z)(x',y',z') \, = \, (x + e^{-z}x' , y + e^zy' , z + z') ,
\end{equation*}
and with Riemannian structure
\begin{equation*}
ds^2 \, = \,  e^{2z}dx^2 + e^{-2z}dy^2 + dz^2 .
\end{equation*}
The fundamental group $\Gamma$ of a Sol-manifold
contains a subgroup of finite index $\Gamma_0$
which fits in a short exact sequence
\begin{equation*}
\{1\} \, \longrightarrow \, 
K  \, \longrightarrow \, 
\Gamma_0  \, \longrightarrow \, 
Q  \, \longrightarrow \,  \{1\} 
\end{equation*}
where $K$ [respectively $Q$]
is a discrete subgroup of the group of isometries
of a Euclidean plane [respectively of  a line];
see Proposition 4.7.7 in \cite{Thur--97}.
In particular, $\Gamma$ has a solvable subgroup of finite index.

\vskip.2cm
 
(\ref{sectionparticularcases}.3)
A \emph{hyperbolic structure} on a manifold $M$ of dimension $n \ge 2$,
here for simplicity assumed to be without boundary,
is a complete Riemannian structure 
of which the sectional curvature is constant of value $-1$.
The universal cover of $M$ can be identified with 
the hyperbolic space $\mathbf H^n$ of dimension $n$.
The fundamental group  $\Gamma$ of $M$ can be identified
with a torsion-free discrete subgroup 
of the group of isometries $\operatorname{Is}(\mathbf H^n)$,
and therefore also with the corresponding group of homeomorphisms
of the boundary $\partial \mathbf H^n \approx \mathbf S^{n-1}$.
\par

A manifold is \emph{hyperbolic} if it has (or more generally can be given)
a hyperbolic structure.
\par

A discrete group of isometries of $\mathbf H^n$ is \emph{elementary} 
if there exists 
either a $\Gamma$-fixed point in $\mathbf H^n \cup \partial \mathbf H^n$,
or a $\Gamma$-invariant pair $\{\alpha, \omega\} \subset \partial \mathbf H^n$;
a discrete group of isometries of $\mathbf H^n$ is elementary
if and only if it has a free abelian subgroup of finite index.
A group $\Gamma = \pi_1(M)$ is never elementary 
if the hyperbolic manifold $M$ is of finite volume, 
and a fortiori if $M$ is compact.

\begin{Lemma}
\label{hyperbolicelementary} 
Let $M$ be an orientable compact $3$-manifold
and let $\Gamma$ denote its fundamental group.
\par
If $M$ is a closed hyperbolic manifold, then $\Gamma$ is not elementary.
\par
If $M$ is a manifold with boundary such that $M \smallsetminus \partial M$
has a hyperbolic structure and if $\Gamma$ is elementary, 
then $M$ is diffeomorphic to one of the following four manifolds:
the ball $\mathbf B^3$, the solid torus $\mathbf S^1 \times \mathbf D^2$,
the thick torus $\mathbf T^2 \times \mathbf I$, 
the non-trivial $\mathbf I$-bundle $\mathbf K^2 \widetilde{\times} \mathbf I$.
\par
In particular, if $\Gamma$ is elementary and $\Gamma \ne \{1\}$, 
then $M$ is a Seifert manifold.
\end{Lemma}

\emph{Proof.}
For general facts on elementary groups, see for example \cite{Ratc--06}.
Suppose that $\Gamma$ is elementary.
Note that $\Gamma$ has no torsion
(=~no elliptic element), 
since $M \smallsetminus \partial M$
is covered by the space $\mathbf H^3$ which is contractible.
\par
If $\Gamma$ contains a hyperbolic element $h_0$,
then $\Gamma$ is infinite cyclic, with $(h_0)^{\mathbf Z}$ of finite index in $\Gamma$
(see also \cite[Theorem 10.7]{Hemp--76}).
It follows that $\Gamma \backslash \mathbf H^3$ 
is diffeomorphic to $\mathbf R \times \mathbf H^2$,
and that $M$ is a solid torus.
\par
If $\Gamma$ contains a parabolic element $p_0$,
there are two possibilities.
Either $\Gamma$ is infinite cyclic and $M$ is again a solid torus.
Or $\Gamma$ has a subgroup of finite index isomorphic to $\mathbf Z^2$;
it follows that $\Gamma \backslash \mathbf H^3$ 
is diffeomorphic to a $\mathbf R$-bundle over a flat Euclidean compact manifold,
and that $M$ is a $\mathbf I$-bundle over one of $\mathbf T^2$, $\mathbf K^2$
(see also \cite[Theorem 10.6]{Hemp--76}).
\par
Otherwise $\Gamma = \{1\}$ and $M$ is a $3$-ball.
\hfill $\square$

\vskip.2cm

(\ref{sectionparticularcases}.4)
A $3$-manifold $M$ is \emph{irreducible}
if any $2$-sphere embedded in $M$ bounds a $3$-ball. 
For example, the sphere $\mathbf S^3$ and the ball $\mathbf B^3$
are irreducible (a theorem of Alexander);
so are the solid torus and the solid Klein bottle.
There are two $\mathbf S^2$-bundles over $\mathbf S^1$:
the trivial bundle $\mathbf S^1 \times \mathbf S^2$, which is orientable,
and the non-trivial bundle, denoted by $\mathbf S^1 \widetilde \times \mathbf S^2$;
both are reducible $3$-manifolds.
Also, non-trivial connected sums (see (\ref{sectiondecomposition}.1) below) 
are reducible.
\par

A surface \emph{properly embedded} in $M$
is a surface $S$ embedded in $M$ such that $\partial S = S \cap \partial M$.
Let $S$ be a surface 
which is either properly embedded  in $M$ and two-sided,
or inside $\partial M$;
then $S$ is  \emph{incompressible} if,
for any simple closed curve $\gamma$ in $S$
and for any disc $D$ embedded in $M$ with $\partial D = D \cap S = \gamma$,
there exists a disc $D'$ embedded in $S$ such that $\partial D' = \gamma$. 
Trivial examples: discs and spheres are always incompressible;
note that our definition of ``incompressible'' is different from 
that in \cite{Hemp--76} in case $S$ is a $2$-disc or a $2$-sphere, 
or from that in \cite[Page 1]{JaSh--79} in case $S$ is a $2$-sphere.
%
\par

Let $S$ be a surface in $M$ which is 
either properly embedded and two-sided or in $\partial M$.
As a consequence of the ``loop + Dehn'' theorem,
$S$ is incompressible if and only if
the homomorphism $\pi_1(S) \longrightarrow \pi_1(M)$
induced by the inclusion is injective
\cite[Chapters 4 and 6]{Hemp--76}. 
\par

A $3$-manifold $M$ is \emph{$\partial$-irreducible} if it is irreducible 
and if the connected components of $\partial M$ (if any) are incompressible.
The exterior of a non-trivial knot 
is $\partial$-irreducible, by Dehn's lemma.
For examples of irreducible manifolds which are not $\partial$-irreducible,
see Remark \ref{RemarkstorusKlein}.

\vskip.2cm

(\ref{sectionparticularcases}.5)
Let $S_0$ be a surface properly embedded in $M$
and let $S_1$ be a surface which is either properly embedded in $M$
or in $\partial M$.
Then $S_0$ and $S_1$ are \emph{parallel}
if there exists an embedding of a thickened surface
$\psi : S \times [0,1] \longrightarrow M$
such that $\psi(S \times \{j\}) = S_j$ (for $j = 0,1$)
and $\psi(\partial S \times [0,1]) \subset \partial M$.
In case $\partial M \ne \emptyset$, 
a surface $S_0$ properly embedded in $M$ is \emph{boundary-parallel}
if there exists a manifold $S_1$ embedded in $\partial M$
such that $S_0$ and $S_1$ are parallel.
For example $\mathbf T^2 \times \{\frac{1}{2}\}$
is boundary parallel in the thickened torus $\mathbf T^2 \times [0,1]$.
\par

An \emph{essential torus} in $M$ is a torus embedded in $M$
which is disjoint from $\partial M$, incompressible, and not boundary parallel.

\vskip.2cm

(\ref{sectionparticularcases}.6)
\emph{Note on irreducibility, boundaries, and compactness.}
Seifert manifolds may have boundaries;
they are compact, by definition.
A Seifert manifold is either irreducible 
or diffeomorphic to one of 
$\mathbf S^1  \times \mathbf S^2$,
$\mathbf S^1 \widetilde \times \mathbf S^2$, 
$\mathbf P^3 \sharp \mathbf P^3$ 
\cite[Proposition 1.12]{Hatc--00}).
An irreducible Seifert manifold is either $\partial$-irreducible
or a solid torus or a solid Klein bottle
\cite[Corollary 3.3]{Scot--83}.
\par

Let $M$ be a complete Riemannian $3$-manifold without boundary.
If $M$ is either Sol or hyperbolic, $M$ is irreducible
\cite[Theorem 2.2]{Bona--02}.
A fortiori, this holds for closed Riemannian $3$-manifolds.
(For Sol and hyperbolic manifolds with boundaries, 
see Section 2.5 in \cite{Bona--02}.)

\begin{Prop}
\label{examples-3-manif}
Let $M$ be a $3$-manifold and let $\Gamma$ denote its fundamental group.
We assume that $\Gamma \ne \{1\}$.
\begin{itemize}
\item[(i)]
If $M$ is a Seifert manifold or if the interior of $M$ is a Sol-manifold, 
then $\Gamma$ is not C$\,{}^*$-simple.
\item[(ii)]
Suppose that the interior of $M$ has a hyperbolic structure;
suppose also that $\Gamma$ is non-elementary
(for elementary groups, see Lemma \ref{hyperbolicelementary}).
Then $\Gamma$ is a strongly Powers group.
(This carries over to any dimension $n \ge 2$.)
\item[(iii)]
Assume that $M$ is compact, orientable, $\partial$-irreducible, 
and contains an essential torus;
assume moreover that $M$ is not a Seifert manifold 
and that the interior of $M$ is not a Sol-manifold.
Then $\Gamma$ is a strongly Powers group.
\end{itemize}
\end{Prop}

\noindent
\emph{Proof of the first two claims.}
(i)
If $M$ is a Seifert manifold, $\Gamma$ either is finite
or contains a normal subgroup isomorphic to $\mathbf Z$;
if the interior of $M$ is a Sol-manifold, 
$\Gamma$ has a solvable subgroup of finite index.
In all cases, $\Gamma$ has an amenable normal subgroup $N \ne \{1\}$,
and therefore cannot be C$\,{}^*$-simple (Item \ref{ItemPowers}.vi).
\par

(ii)
Let $\Gamma$ be as in (ii), viewed as a group of isometries
of the hyperbolic space $\mathbf H^n$.
We can apply Proposition \ref{Powersgeneral}, 
with $\Omega = \partial \mathbf H^n \approx \mathbf S^{n-1}$.
The  set $\overline{L_{\Gamma}}$ as defined before Proposition \ref{Powersgeneral}
coincides  with the usual ``limit set'' of $\Gamma$.
It is easy to check that the action of $\Gamma$ is strongly hyperbolic on $\mathbf S^{n-1}$
and strongly faithful on $\overline{L_{\Gamma}}$.
To apply Proposition \ref{Powersgeneral},
let us check that any normal subgroup $N \ne \{1\}$ of $\Gamma$
contains hyperbolic elements.
\par

On the one hand, $N$ cannot contain non-identity elliptic transformations,
because $\Gamma$ is torsion-free 
and discrete in $\operatorname{Is}(\mathbf H^n)$.
On the other hand,
if $N$ contains parabolic transformations,
then the closure of
\begin{equation*} 
\{ \eta \in \mathbf S^{n-1}
\hskip.1cm \vert \hskip.1cm
\gamma \eta = \eta 
\hskip.2cm \text{for some parabolic} \hskip.2cm
\gamma \in N \} 
\end{equation*}
coincides with $\overline{L_{\Gamma}}$,  by minimality.
Hence $N$ contains a pair $(p_1, p_2)$ of parabolic elements
with distinct fixed points.
Then $p_1^k p_2^k$ is hyperbolic for $k$ large enough.
Thus, in all cases,  $N$ contains hyperbolic elements.
\hfill $\square$

\vskip.2cm

Similarly, in a non-elementary torsion-free Gromov-hyperbolic group $\Gamma$,
any subnormal subgroup which is not $\{1\}$
contains hyperbolic elements,
so that $\Gamma$ is a strongly Powers group.

\vskip.2cm

We will complete the proof of Claim (iii) of Proposition \ref{examples-3-manif} 
with Lemma \ref{longlemme}. 
Before this, we need some preliminaries.

\begin{Lemma}
\label{L1duchap6}
Let $N$ be a compact orientable $3$-manifold.
Assume that $N$ is $\partial$-irreducible 
and that there exists a connected component $T$ of its boundary
which is a $2$-torus.
Assume moreover that the fundamental group $\pi_1(N)$
contains a normal subgroup $K \ne \{1\}$
which is peripheral, and more precisely in the image of $\pi_1(T)$.
\par
Then one of the following conclusions holds:
\begin{itemize}
\item[(i)]
$K \approx \mathbf Z$ and $N$ is a Seifert manifold,
\item[(ii)]
$K \approx \mathbf Z^2$ and $N \approx \mathbf T^2 \times \mathbf I$,
\item[(iii)]
$K \approx \mathbf Z^2$ and $N \approx \mathbf K^2 \widetilde{\times} \mathbf I$.
\end{itemize}
\end{Lemma}

\noindent
\emph{Proof.}
Note that $\pi_1(T) \approx \mathbf Z^2$ can be identified
with a subgroup of $\pi_1(N)$, 
because $T$ is by assumption incompressible in $N$.
Observe that $K$, which is a subgroup of $\pi_1(T)$,
is isomorphic to one of $\mathbf Z$ or $\mathbf Z^2$.
Consider two cases,  depending on the size of $Q \Doteq \pi_1(N)/K$.
\par

\emph{Case 1:  $Q$ is finite.}
Since $K$ is a group of finite index
in the group $\pi_1(N)$ containing
$\pi_1(T) \approx \mathbf Z^2$,
we cannot have $K \approx \mathbf Z$;
hence $K \approx \mathbf Z^2$.
Theorem 10.6 in \cite{Hemp--76} implies
that $N$ is a $\mathbf I$-bundle over a surface, say $F$.
Since $\pi_1(F)$ has $\mathbf Z^2$ as a subgroup of finite index,
$F$ is either a $2$-torus or a Klein bottle.
Hence the situation is as in (ii) or (iii).
\par

\emph{Case 2: $Q$ is infinite.}
By Consequence (3) in Theorem 11.1 of \cite{Hemp--76},
we have $K \approx \mathbf Z$.
By the theorem of Seifert fibration
(see for our case Theorem II.6.4 in \cite{JaSh--79}),
$N$ is a Seifert manifold.\hfill $\square$

\vskip.2cm

Let $M$ be a compact $3$-manifold 
which contains a two-sided properly embedded surface $S$
(this $S$ need not be connected).
\begin{itemize}
\item[$\bullet$]
The manifold $M_S^*$ obtained by \emph{splitting $M$ along $S$}
is the complement in $M$ of a regular open neighbourhood of $S$.
\end{itemize}
(Our $M_S^*$ is that manifold which is denoted
by $\sigma_S(M)$ in \cite[Page 4]{JaSh--79}.)
If $M$ is orientable,
$M_S^*$ is irreducible if and only if $M$ is irreducible;
this is a theorem of Waldhausen \cite[Theorem 1.8]{Wa--67ab}.
It follows that, if $M$ is $\partial$-irreducible and if $S$ is incompressible, 
then  $M_S^*$ is $\partial$-irreducible.

\begin{Lemma}
\label{longlemme}
Let $M$ be a manifold which is compact, orientable, $\partial$-irreducible,
and which contains an essential torus $T$.
Assume moreover that $M$ is not a Seifert manifold, and the interior of $M$
is not a Sol-manifold.

\par

Then the fundamental group $\Gamma$ of $\pi_1(M)$
is a strongly Powers group.
\end{Lemma}

\noindent
\emph{Proof of Lemma \ref{longlemme} in case $T$ is separating.}
Denote by $M_1, M_2$ 
the two components of the manifold $M_T^*$ 
obtained by splitting $M$ along $T$.
Since $M$ is $\partial$-irreducible and $T$ is essential,
$M_1$ and $M_2$ are also $\partial$-irreducible.
Set $A = \pi_1(M_1)$, $B = \pi_1(M_2)$, 
and $C = \pi_1(T) \approx \mathbf Z^2$,
so that $\Gamma = A \ast_C B$
by the Seifert--van Kampen theorem
(for which the reference we like best is $\S$ 3 of \cite{Rham--69});
observe that the homomorphism $C \longrightarrow A$
induced by the inclusion $T \subset M_1$ is injective,
because $T$ is incompressible;
the same holds for $C \longrightarrow B$.

The amalgam is non-trivial.
Indeed, if $C = A$, then $\pi_1(M_1) \approx \mathbf Z^2$, 
and  Lemma \ref{L1duchap6} 
shows that $M_1$ is one of
$\mathbf T^2 \times \mathbf I$,  $\mathbf K^2 \widetilde{\times} \mathbf I$.
But the first case is ruled out by ``$T$ essential''
and the second case is ruled out because
$\pi_1(\mathbf K^2 \widetilde{\times} \mathbf I) = \pi_1(\mathbf K^2)
\not\approx \mathbf Z^2$.

We claim that the amalgam is non-degenerate.
For suppose the opposite, namely that $[A:C] = [B:C] = 2$.
Case (ii) if Lemma \ref{L1duchap6} 
is ruled out, so that $M_1$ and $M_2$ are diffeomorphic to
$ \mathbf K^2 \widetilde{\times} \mathbf I$.
Then $\partial M_1 \approx \partial M_2 \approx \mathbf T^2$,
and we have two $2$-coverings $N_j \longrightarrow M_j$
with $N_j \approx \mathbf T^2 \times \mathbf I$ ($j=1,2$).
Thus we have a $2$-covering $N \longrightarrow M$
where $N = N_1 \cup N_2$ is a $\mathbf T^2$-bundle over $\mathbf S^1$,
so that $N$ can be viewed as obtained from
$\mathbf T^2 \times \mathbf I$ by identifying 
the two connected components of the boundary
with a gluing diffeomorphism  of $\mathbf T^2$.
A torus bundle over $\mathbf S^1$ is 
Sol if the gluing map is Anosov,
and Seifert if the gluing map is periodic or reducible
(Theorem 5.5 in \cite{Scot--83} 
or Exercice 3.8.10 in \cite{Thur--97}). 
If $N$ is Sol, then  $M$ is also Sol
\cite[Theorem 5.3.i]{Scot--83}.
If  $N$ is Seifert, then  $M$ is also Seifert
\cite[Theorem II.6.3]{JaSh--79}.
Thus, in all cases, $M$ is Seifert or Sol,
but this has been ruled out.
Thus this case does not occur.
\par

By Claim (i) of Theorem~\ref{MainResult},  
either $\Gamma$ is a strongly Powers group,
or $C_k \ne \{1\}$ for all $k \ge 1$.
From now on, we will assume that $C_2 \ne \{1\}$
(so that $C_1 \ne \{1\}$),
and we will obtain a contradiction.
\par

Set $C_A = \bigcap_{a \in A} a^{-1} C a$ 
and $C_B = \bigcap_{b \in B} b^{-1} Cb$,
so that $C_1 = C_A \cap C_B$
is a non-trivial subgroup of $C \approx \mathbf Z^2$.
Observe that $C_A$ is isomorphic to 
either $\mathbf Z$ or $\mathbf Z^2$, 
and also that $C_A$ is the largest subgroup of $C$
which is normal in $A$.
\par

Note that neither $M_1$ nor $M_2$ can be
diffeomorphic to $\mathbf T^2 \times \mathbf I$,
otherwise $T$ would be boundary parallel in $M$,
and this is not the case by hypothesis.
Hence each of $M_1, M_2$
is either as in (i) or as in (iii)
of Lemma \ref{L1duchap6}, 
and we will end the proof by showing that each case is ruled out.
\par

We have already seen that the case
$M_1 \approx M_2 \approx \mathbf K^2 \widetilde{\times} \mathbf I$
cannot occur.

Suppose then that $M_1$ and $M_2$ are Seifert,
with 
\begin{equation*}
A \,  \rhd \, C_A \, \approx \, \mathbf Z \, \approx \, C_B \, \lhd \, B .
\end{equation*}
Since any subgroup of $\mathbf Z$ is characteristic,
$C_1 = C_A \cap C_B$ is normal in both $A$ and $B$,
so that $C_1 \approx \mathbf Z$ is also normal in $\Gamma$.
By the theorem of Seifert fibrations 
(already used in the proof of Lemma \ref{L1duchap6}), 
it follows that $M$ is Seifert, but this has been ruled out.
\par

We can finally suppose that $M_1$ is Seifert
with $Z \approx C_A \lhd A$,
and $M_2 = \mathbf K^2 \widetilde{\times} \mathbf I$
with $\partial M_2 \approx \mathbf T^2$
and $C = \pi_1(\partial M_2) \approx \mathbf Z^2$;
since $C_B$ is the maximal subgroup of $C$
which is normal in $B$, we have $C_B = C$.
It follows that $C_1 = C_A \cap C_B = C_A$.
Since $C_1$ is normal in $A$,
\begin{equation*}
\{1\} \, \ne \, C_2 \, = \, 
C_1 \cap \Big( \bigcap_{b \in B} b^{-1} C_1 b \Big) \, \lhd \,
C_A \, \approx \, \mathbf Z ,
\end{equation*}
$C_2$ is normal in both $A$ and $B$,
and therefore also in $\Gamma$.
(Though we will not use this, let us observe that
$C_k = C_2$ for all $k \ge 2$.)
It follows as above that $M$ is Seifert,
which has been ruled out.
\hfill $\square$

\vskip.2cm

\noindent
\emph{Proof of Lemma \ref{longlemme} in case $T$ is non-separating.}
Consider again $M_T^*$.
Since $M$ is $\partial$-irreducible and $T$ is essential,
$M_T^*$ is also $\partial$-irreducible.
Set $G = \pi_1(M_T^*)$;
the two boundary components of $M_T^*$ coming from the splitting
correspond to two isomorphic subgroups of $G$, 
say $H$ and $\theta(H)$ where $\theta$ 
is an isomorphism with domain $H$,
and $H \approx \pi_1(T) \approx \mathbf Z^2$.
If one had $H = G$ or $\theta (H) = G$,
the manifold $M$ would be a $\mathbf T^2$-bundle 
over $\mathbf S^1$ \cite[Theorem 10.2]{Hemp--76}, 
and this cannot be (as above).
It follows that the HNN-extension
$\Gamma = \pi_1(M) = HNN(G,H,\theta)$
is non-degenerate.
By  Claim (ii) of Theorem~\ref{MainResult},
either $\Gamma$ is a strongly Powers group,
or $H_k \ne \{1\}$ for all $k \ge 1$.
From now on, we will assume that $H_2 \ne \{1\}$
(so that $H_1 \ne \{1\}$),
we will obtain a contradiction.
\par

Set $H_G = \bigcap_{g \in G} g^{-1} H g$,
so that $H_1 = H_G \cap \tau^{-1} H \tau \cap \tau H \tau^{-1}$
is a non-trivial subgroup of $H \approx \mathbf Z^2$.
Observe that $H_G$ is isomorphic to
either $\mathbf Z$ or $\mathbf Z^2$,
and also that $H_G$ is the maximal subgroup of $H$
which is normal in $G$.
\par

The boundary of $M_T^*$ has at least 
two connected components which are $2$-tori. Since 
$\partial(\mathbf K^2 \widetilde{\times} \mathbf I) \approx \mathbf T^2$
is connected, $M_T^*$ cannot fit in Case (iii) of Lemma \ref{L1duchap6}.
If one had $H_G \approx \mathbf Z^2$ 
and $M_T^* \approx \mathbf T^2 \times \mathbf I$, 
the manifold $M$ would be a $\mathbf T^2$-bundle over $\mathbf S^1$,
consequently would be Seifert or Sol,
and this is ruled out (as above).
\par

Therefore, because of Lemma \ref{L1duchap6},
we can assume that $M_T^*$ is Seifert and that $H_G \approx \mathbf Z$.
Choose a generator $h$ of $H_1$, so that $H_1 = \langle h \rangle$.
Since $\langle h \rangle$ is characteristic in $H_G$,
the subgroup $\langle h \rangle$ is normal in $G$.
Since
\begin{equation*}
\{1\} \, \ne \, H_2 \, = \, 
\langle h \rangle \cap \tau^{-1} \langle h \rangle \tau \cap
\tau \langle h \rangle \tau^{-1} ,
\end{equation*}
there exists a pair of non-zero integers $p,q$ such that
$\tau^{-1} h^p \tau = h^q$.
By what we know about Baumslag-Solitar subgroups
of $3$-manifold groups \cite[Theorem VI.2.1]{JaSh--79},
this implies that $q = \pm p$.
Hence $\Gamma$ contains a normal subgroup
$\langle h^p \rangle \approx \mathbf Z$.
It follows as above that $M$ is a Seifert manifold,
and this has been ruled out.
\par

This ends the proof of Lemma \ref{longlemme},
and thus  of Proposition \ref{examples-3-manif}.
\hfill $\square$

\vskip.2cm

We will need the following lemma in Section \ref{sectiondecomposition}.
Its proof makes use of the notion of splitting recalled just before
Lemma \ref{longlemme}. 

\begin{Lemma}
\label{torusKlein}
Let $M$ be a $3$-manifold with $\pi_1(M) \ne \{1\}$.
Assume that $M$ is irreducible and not $\partial$-irreducible.
\par
Then
either $\pi_1(M)$ is a non-trivial free product,
or $M$ is a solid torus or a solid Klein bottle.
\end{Lemma}

\noindent 
\emph{Proof.}~Since $\partial M$ has a compressible component, say $F$,
there exists a $2$-disc $D$ in $M$ such that
$\partial D = D \cap F$,
and $\partial D$ not contractible in~$F$.
Let $M_D^*$ be the result of splitting $M$ along $D$.
There are two cases to distinguish.
\par

Suppose first that  $M_D^*$ is connected.
By the Seifert--van Kampen theorem,
we have $\pi_1(M) = \pi_1(M_D^*) \ast \mathbf Z$.
If $M_D^*$ is simply connected,
then $M_D^*$ is diffeomorphic to a $3$-ball
(otherwise $M$ would not be irreducible),
and $M$ is diffeomorphic to either a solid torus
or a solid Klein bottle, depending on the action
of the gluing map on the orientation of $D$.
If $M_D^*$ is not simply connected,
then $\pi_1(M)$ is a non-trivial free product.
\par

If $M_D^*$ has two connected components,
say $M'$ and $M''$, then $\pi_1(M) = \pi_1(M') \ast \pi_1(M'')$,
again by the Seifert--van Kampen theorem.
Observe that neither $M'$ nor $M''$ is simply connected
(otherwise one of $M'$, $M''$ would be a $3$-ball,
and $\partial D$ would be contractible in $F$),
so that $\pi_1(M)$ is a non-trivial free product.
\hfill $\square$

\begin{Rem}
\label{RemarkstorusKlein}
Each case in the conclusion of Lemma \ref{torusKlein} 
occurs.
\end{Rem}

\noindent \emph{Proof.}
Let $M_1$ be a compact $3$-manifold 
which is irreducible and which has a non-empty boundary.
Let $D,D'$ be two disjoint $2$-discs in $\partial M_1$,
and let $M$ be the result of attaching a handle
$\mathbf D^2 \times \mathbf I$ to $D$ and $D'$.
Then $\pi_1(M) \approx \pi_1(M_1) \ast \mathbf Z$.
The manifold $M$ is irreducible and, 
if  $M_1$ is not simply connected, 
$\pi_1(M)$ is a non-trivial free product; 
observe that $M$ is not $\partial$-irreducible.
\hfill $\square$

\section{\textbf{Reduction to Proposition \ref{examples-3-manif} 
}}
\label{sectiondecomposition}

(\ref{sectiondecomposition}.1)
Our first reduction step is inspired by 
the Kneser-Milnor theorem on decompositions of $3$-manifolds by connected sums.
For simplicity, we assume from now on that 
\begin{center}
all $3$-manifolds below are \emph{orientable}
\end{center}
(with the exception of Proposition \ref{nonorientable}). 
We begin with some reminder.
\par

Given two $3$-manifolds $M_1, M_2$, one can define 
the \emph{connected sum} $M = M_1 \sharp M_2$,
with $\pi_1(M)$ the free product $\pi_1(M_1) \ast \pi_1(M_2)$.
\begin{itemize}
\item[$\bullet$]
A $3$-manifold $M$ is \emph{prime} if, 
whenever $M$ is diffeomorphic to a connected sum $M_1 \sharp M_2$,
at least one of $M_1, M_2$ is diffeomorphic to $\mathbf S^3$.
\end{itemize}
It is a standard result that a compact manifold which is prime
is either  a $\mathbf S^2$-bundle over $\mathbf S^1$
(see \ref{sectionparticularcases}.4)
or irreducible \cite[Lemma 3.13]{Hemp--76}.

Any compact $3$-manifold $M$ has  a decomposition in connected sum
\begin{equation}
\label{eqconnectedsumM}
M \, \approx \, M_1 \sharp \cdots \sharp M_k,
\end{equation}
where each $M_j$ is prime,
and we have a free product decomposition
\begin{equation}
\label{eqconnectedsumG}
\pi_1(M) \, \approx \,  \pi_1(M_1) \ast \cdots \ast \pi_1(M_k) .
\end{equation}
This is the \emph{Kneser-Milnor decomposition} of~$M$
\cite[Chapter 3]{Hemp--76}. 
Unless $M$ itself is a $3$-sphere, 
we ask moreover that no $M_j$ is a $3$-sphere.
\par

Let us discuss the possible factors with $\pi_1(M_j) = \{1\}$ in (\ref{eqconnectedsumG}).
Standard arguments of algebraic topology show that this can occur only if
$M_j$ is a $3$-ball or a homotopy $3$-sphere.
On the one hand, to avoid $3$-balls,
we consider as in Section \ref{section1}
the manifold $\widehat M$ obtained from $M$ by filling 
spherical boundary components of $\partial M$ with $3$-balls
\cite[Lemma 3.7]{Hemp--76};
note that  $\pi_1(\widehat M) \approx \pi_1(M)$
by the Seifert--van Kampen theorem.
On the other hand, we know that no exotic homotopy sphere can occur,
by Perelman's proof of the Poincar\'e conjecture.
(Thus our $\widehat M$ is  identical to the ``Poincar\'e completion''
$\mathcal P (M)$ of \cite{Hemp--76}.)
\par

\begin{Prop}[\textbf{first reduction}]
\label{first reduction}
Let $M$ be a compact orientable $3$-manifold 
and let $\Gamma$ be its fundamental group.
Then one of the following statements is true:
\begin{itemize}
\item[(i)] 
$\Gamma$ is a non-degenerate free product,
and is therefore a strongly Powers group;
\item[(ii)]
Either $\widehat M$ is a Seifert manifold, 
or the interior of $M$ is a Sol-manifold;
and therefore $\Gamma$ is not C$\, {}^*$-simple;
\item[(iii)] 
$\hat M$ is $\partial$-irreducible, and is not a Seifert manifold,
and the interior of $M$ is not a Sol-manifold.
\end{itemize}
\end{Prop}

\noindent
\emph{Proof.}
Suppose first that $\widehat M$ is not prime.
Then $\Gamma$ is a non-trival free product, as discussed above,
so that $\Gamma$ is either as in (i),
or an infinite dihedral group.
In the second case, it is a result
of \cite{Tao--62}
that $\widehat M$ is a connected sum $\mathbf P^3 \sharp \mathbf P^3$
of two projective spaces; this is a circle bundle over $\mathbf P^2$,
and in particular a Seifert manifold as in (ii).
\par

Suppose now that $\widehat M$ is prime and reducible.
Then $\widehat  M \approx \mathbf S^1 \times \mathbf S^2$
is a Seifert manifold.
\par

Suppose finally that $\widehat M$ is irreducible.
If $\pi_1(M) = \{1\}$,  then $\widehat M \approx \mathbf S^3$
is a Seifert manifold 
(recall that, for simplicity of our statement and as in (\ref{sectionparticularcases}.1),
we rely on Perelman's work).
If $\pi_1(M) \ne \{1\}$, then $\widehat  M$ is 
 either a Seifert manifold or $\partial$-irreducible,
by Lemma \ref{torusKlein}. 
\hfill $\square$

\vskip.2cm

\emph{Remark:}
the three cases of the proposition are exclusive,
because, if $\pi_1(M)$ is a non-trivial free product, 
then $\widehat M$ is a non-trivial connected sum
(Stalling's proof of Kneser's conjecture \cite[Chapter 7]{Hemp--76}).

%

\vskip.2cm

(\ref{sectiondecomposition}.2)
Our second reduction step relies
on the JSJ decomposition of $3$-manifolds by essential tori,
that we recall below.
An irreducible manifold $M$ is \emph{atoroidal}
if any incompressible torus in $M$ 
is parallel to a component of $\partial M$.
Here is the basic \emph{JSJ decomposition} theorem 
of Jaco-Shalen and Johannson, announced by Waldhausen
(see \cite{Wald--69},  \cite[Page 157]{JaSh--79},
\cite[Page 483]{Scot--83}, and \cite[Theorem 3.4]{Bona--02}):

\begin{Thm}[\textbf{JSJ}]
\label{JSJ}
Let $M$ be a compact orientable $3$-manifold
which is  $\partial$-irreducible.
\par
There exists a minimal finite family of disjoint two-sided
essential tori $T_1, \hdots, T_k$ such that 
each connected component of 
$M \smallsetminus \bigcup_{j=1}^k T_j$
is either atoroidal or a Seifert manifold,
and this family is unique up to isotopy.
\end{Thm}

The components of $M \smallsetminus \bigcup_{j=1}^k T_j$
are called the \emph{pieces}.
Using the same notation, we have the following corollary.

\begin{Cor}
There is a graph of groups $\mathbf G  = (G,Y)$,
of which the graph $Y$ has one vertex for each piece 
and one geometric edge for each torus,
with $G_e = \pi_1(\mathbf T^2) \approx \mathbf Z^2$
for all $e \in E(Y)$,
such that $\pi_1(M)$ is isomorphic 
to the fundamental group of $\mathbf G$.
\par
Moreover, if $k > 0$, 
either $M$ is a torus bundle over $\mathbf S^1$ (and then $k=1$),
or, for each edge $e \in E(Y)$,
the images of $G_e$ by $\varphi_e$ and $\varphi_{\overline{e}}$
are proper subgroups of $G_{t(e)}$ and $G_{s(e)}$ respectively.
\end{Cor}

\noindent \emph{Proof.}
For the second part of the corollary, 
suppose that there is an edge $e$ in $Y$ 
with source $s$ and terminus $t$,
corresponding respectively to a torus $T_e$ and to pieces $V_s$ and $V_t$, 
such that $\varphi_e : \pi_1(T_e) \longrightarrow \pi_1(T_t)$ is an isomorphism.
Then $V_t$ is a thickened torus $T_e \times \mathbf I$
by \cite[Theorem 10.2]{Hemp--76}.
Since two $T_j$ are never parallel, 
this implies $s=t$ and $k=1$,
so that $M$ is a torus bundle over $\mathbf S^1$.
\hfill $\square$

\vskip.2cm

The following proposition is a straightforward consequence of Theorem \ref{JSJ}. 

\begin{Prop}[\textbf{second reduction}]
\label{second reduction}
Let $M$ be a compact orientable $3$-manifold
which is  $\partial$-irreducible.
Then one of the following statements is true:

(i)
$M$ is an atoroidal manifold.

(ii)
$M$ is a Seifert manifold or the interior of $M$ is a Sol-manifold.

(iii)
$M$ contains a two-sided essential embedded torus,
$M$ is not a Seifert manifold 
and the interior of $M$ is not a Sol-manifold.
(In other words $M$ fulfills the hypothesis of 
Claim (iii) in Proposition \ref{examples-3-manif}.) 
\end{Prop}

(\ref{sectiondecomposition}.3)
In three preprints made public in 2002 and 2003, 
Grisha Perelman, following the Hamilton program and
using the Ricci flow, has sketched a proof 
of  Thurston's geometrisation conjecture;
particular cases have been known before \cite{Thur--82},
and details have been provided by other authors.
See \cite{B$^3$MP--10}, \cite{MoTi}, and references there. 
We can state the theorem as follows:

\begin{Thm}[\textbf{Thurston-Perelman}]
\label{TP}
Let $M$ be an irreducible compact orientable $3$-manifold.
If $M$ is atoroidal, then the interior of $M$ has a hyperbolic structure.
\end{Thm} 

\begin{Cor}
\label{CorTP}
Let $M$ be a compact orientable $3$-manifold, with fundamental group $\Gamma$.
Then:
\par
(i) either $\widehat M$ is a Seifert, 
or the interior of $M$ is Sol,
and then $\Gamma$ is not C$\,{}^*$-simple,
\par
(ii) or $\Gamma$ is a strongly Powers group.
\end{Cor}

\noindent \emph{Proof.}
This is a consequence of Propositions \ref{first reduction} 
and \ref{second reduction},  
Theorem \ref{TP}, Corollary \ref{CorTP}, 
Lemma \ref{hyperbolicelementary}, 
and Proposition \ref{examples-3-manif}. 
\hfill $\square$

\vskip.2cm

\noindent
\textbf{End of proof of Claim (iv) in Theorem \ref{MainResult}.} 
The first part of this Claim (iv) is precisely the previous corollary.
For the second part, 
let $K$ be a knot in $\mathbf S^3$.
Denote by $M$ the complement of an open tubular neighbourhood of $K$,
and let $\Gamma = \pi_1(M)$ be the group of $K$.
Then $M$ is irreducible, by the Alexander-Sch\"onflies theorem.
Assume moreover that $K$ is not trivial, 
so that $\partial M$ is an incompressible torus 
by the ``loop + Dehn'' theorem \cite[Chapter 4]{Hemp--76},
and $M$ a $\partial$-irreducible manifold.
\par

Suppose first that the JSJ decomposition of $M$ is trivial,
so that $M$ is either Seifert or atoroidal,
by Proposition \ref{second reduction}.
If $M$ is a Seifert manifold, then 
$\Gamma$ is not C$\, ^*$-simple
and $K$ is a \emph{torus knot};
see \cite{Budn--06}, Proposition 4, first case\footnote{In
the second case of Proposition 4 of \cite{Budn--06},
$n \ge 1$ should be $n \ge 2$.}, 
with $n=1$.
If $M$ is atoroidal, Thurston's hyperbolisation theorem implies
that the interior of $M$ is a hyperbolic manifold of finite volume
(see e.g. \cite{Bona--02}, Section 6.1); 
in this case, $K$ is a \emph{hyperbolic knot} and
$\Gamma$ is a strongly Powers group,
by Proposition \ref{examples-3-manif}.ii. \par 

Suppose now that the JSJ decomposition of $M$ is non-trivial
(i.e. involves at least one torus),
so that $K$ is a \emph{satellite knot}.
Observe that $M$ is not a Seifert manifold
(this would imply that $K$ is a torus knot \cite[Note added in proof]{BuZi--66})
and that the interior of $M$ is not a Sol-manifold
(the list of Sol-manifolds with boundaries is very short,
and does not contain any knot complement \cite[Theorem 2.15]{Bona--02}).
Then $\Gamma$ is a strongly Powers group,
by Lemma \ref{longlemme}. 
\hfill $\square$

\vskip.2cm

For non-orientable manifolds, we will restrict ourselves
to the following simple statement:

\begin{Prop}
\label{nonorientable} 
Let $M$ be a non-orientable connected compact  $3$-manifold, 
$\Gamma$  its fundamental group,
and $\Gamma'$  the fundamental group 
of the total space of the orientation cover of $M$;
in particular,  $\Gamma'$ is a subgroup of index $2$ in $\Gamma$.
\par

Then $\Gamma$ is C$\,{}^*$-simple if and only if 
$\Gamma'$ is C$\,{}^*$-simple.
\end{Prop}

\noindent \emph{Proof.} 
If $\Gamma$ is of order two,
neither $\Gamma$ nor $\Gamma' = \{1\}$ is C$^*$-simple,
and the proposition holds for a trivial reason.
We assume from now on that $\Gamma$ is not of order two,
and it follows that $\Gamma$ is infinite
(a result of D. Epstein,  Theorem 9.5 in \cite{Hemp--76}).

Any subgroup of finite index in a C$^*$-simple group 
is itself C$^*$-simple \cite{BeHa--00}.
In particular, if $\Gamma$ is C$^*$-simple, so is $\Gamma'$.

Assume finally that $\Gamma'$ is C$^*$-simple.
Then $\Gamma'$ is icc  \cite[Appendix J]{Harp--07}.
It follows from Lemma 9.1 and Proposition 3.2 in \cite{HaPr--07} 
that $\Gamma$ is icc,
and then from \cite{BeHa--00} that $\Gamma$ is C$^*$-simple.
\hfill $\square$

\end{document}